\documentclass[11pt, a4paper]{article}
\usepackage{amsmath}
\usepackage{amsfonts}
\usepackage{amssymb}
\usepackage{ulem}

\usepackage{amsthm}
\usepackage{url}

\usepackage{cite}
\usepackage{fancyhdr}

\usepackage{graphicx, color}
\usepackage{hyperref}

\begin{document}

\title{Convex Pentagon Tilings and Heptiamonds, II}
\author{ Teruhisa SUGIMOTO$^{ 1), 2)}$ and Yoshiaki ARAKI$^{2)}$ }
\date{}
\maketitle

{\footnotesize

\begin{center}
$^{1)}$ The Interdisciplinary Institute of Science, Technology and Art

$^{2)}$ Japan Tessellation Design Association

Sugimoto E-mail: ismsugi@gmail.com, Araki E-mail: 
yoshiaki.araki@tessellation.jp
\end{center}

}

\medskip

{\small
\begin{abstract}
\noindent
In the previous manuscript, new tilings (tessellations) were presented 
using convex pentagonal tiles belonging to Type 1 and Type 5. 
The convex pentagon tilings are related to heptiamond tilings. Later, the 
authors found Johannes Hindriks' site that summarized the research results 
of heptiamond tilings (tessellations). In this manuscript, the results of 
converting the heptiamond tilings of Hindriks to convex pentagon tilings are 
introduced. As a result, many new convex pentagon tilings are presented.
\end{abstract}

}

\section{Introduction}
\label{section1}

\begin{figure}[htbp]
 \centering\includegraphics[width=14.5cm,clip]{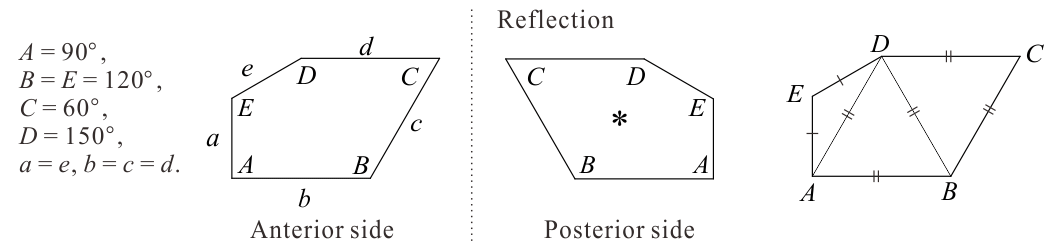} 
  \caption{{\small 
TH-pentagon.} 
\label{fig1}
}
\end{figure}

\begin{figure}[htbp]
 \centering\includegraphics[width=14cm,clip]{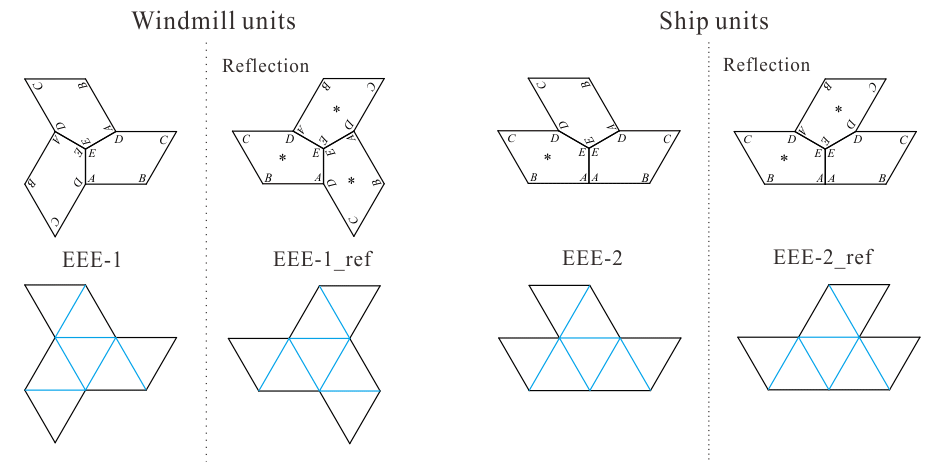} 
  \caption{{\small 
Windmill units and ship units.} 
\label{fig2}
}
\end{figure}

In the previous manuscript \cite{S_and_A_2017} the convex pentagon in 
Figure~\ref{fig1} is called a TH-pentagon, and the tilings (tessellations) 
by TH-pentagons are presented.  When the TH-pentagon generates a 
tiling, the vertex $E$ belongs to 
$A + D + E = 360^ \circ$ or $3E = 360^ \circ$. Furthermore, 
$A + D + E = 360^ \circ$ and $3E = 360^ \circ$ can not coexist within 
the tiling. The tilings with $A + D + E = 360^ \circ$ are always variations 
of Type 1 tilings in Figure 2 in \cite{S_and_A_2017}. The tilings with 
$3E = 360^ \circ $ are formed by windmill units and ship units in Figure~\ref{fig2}. 

Since the windmill unit and the ship unit can be considered as two types 
among 24 heptiamonds, tilings with $3E = 360^ \circ$ by TH-pentagons are 
equivalent to tilings by the two types of heptiamonds~\cite{G_and_S_1987, 
S_and_A_2017}.

\renewcommand{\figurename}{{\small Figure.}}

\renewcommand{\figurename}{{\small Figure.}}

In the previous manuscript \cite{S_and_A_2017}, many convex pentagon tilings by 
windmill units and ship units were presented, but it is probable that there are still 
undiscovered pentagon tilings from the properties of TH-pentagons. This was 
confirmed by Johannes Hindriks' site~\cite{Hindriks-site} that summarizes the research 
results of heptiamond tilings (tessellations). On Hindriks' site~\cite{Hindriks-site}, there 
are many tilings (tessellations) that can be formed by each of 24 types of 
heptiamonds, and the windmill unit and ship unit are Heptiamond 24 (Piece 
24) and Heptiamond 15 (Piece 15), respectively.

Therefore, in this manuscript, tilings of Heptiamond 24 and Heptiamond 15 in 
Hindriks' site~\cite{Hindriks-site} are converted to convex pentagon tilings, 
and their properties are introduced.

\section{Tilings by only the windmill units}
\label{section2}

In this section, the tiling of Hindriks' site~\cite{Hindriks-site} are introduced 
which correspond to tiling by only the windmill unit. The articles related to 
Heptiamond 24 on Hindriks' site~\cite{Hindriks-site} are as follows.

\medskip

37sec824 $\colon$ \url{http://www.jhhindriks.info/37/37sec824.htm}

37sec800 $\colon$ Sections 8.2, 8.3, and 8.5 in 
\url{http://www.jhhindriks.info/37/37sec800.htm}

\medskip

Hindriks calls the figure which corresponds to the hexagonal flower L1 unit 
(HFL1-unit) in \cite{S_and_A_2017} a Miniflower. Hindriks' definition of the 
term \cite{Hindriks-site} is at \url{http://www.jhhindriks.info/37/37terms.htm}. 
Hindriks' site~\cite{Hindriks-site} introduces many tilings with the basic 
shape of Heptiamond 24 called Chevy.

In conclusion, if converting the tilings of Hindriks' Heptiamond 24 to 
convex pentagon tilings, they are all tilings that can be formed according 
to the properties shown in Section 3 in \cite{S_and_A_2017}. That is, a new 
tiling (new filling rule) as the TH-pentagon using only the windmill units does not 
exist on Hindriks' site~\cite{Hindriks-site}.

\begin{figure}[htbp]
 \centering\includegraphics[width=14.5cm,clip]{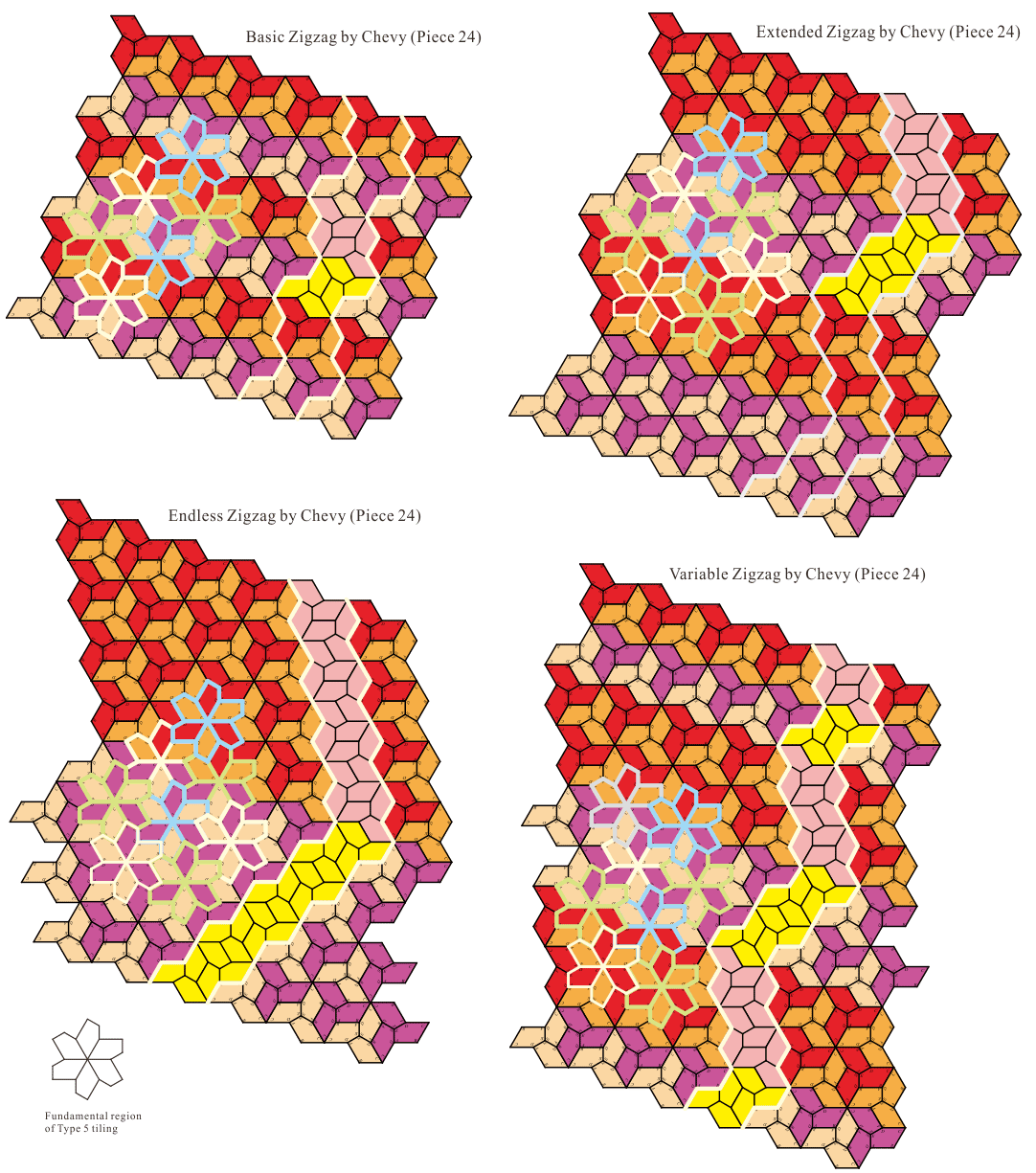} 
  \caption{{\small 
Examples of Zigzag that converted to convex pentagon tiling.} 
\label{fig3}
}
\end{figure}

\begin{figure}[htbp]
 \centering\includegraphics[width=15cm,clip]{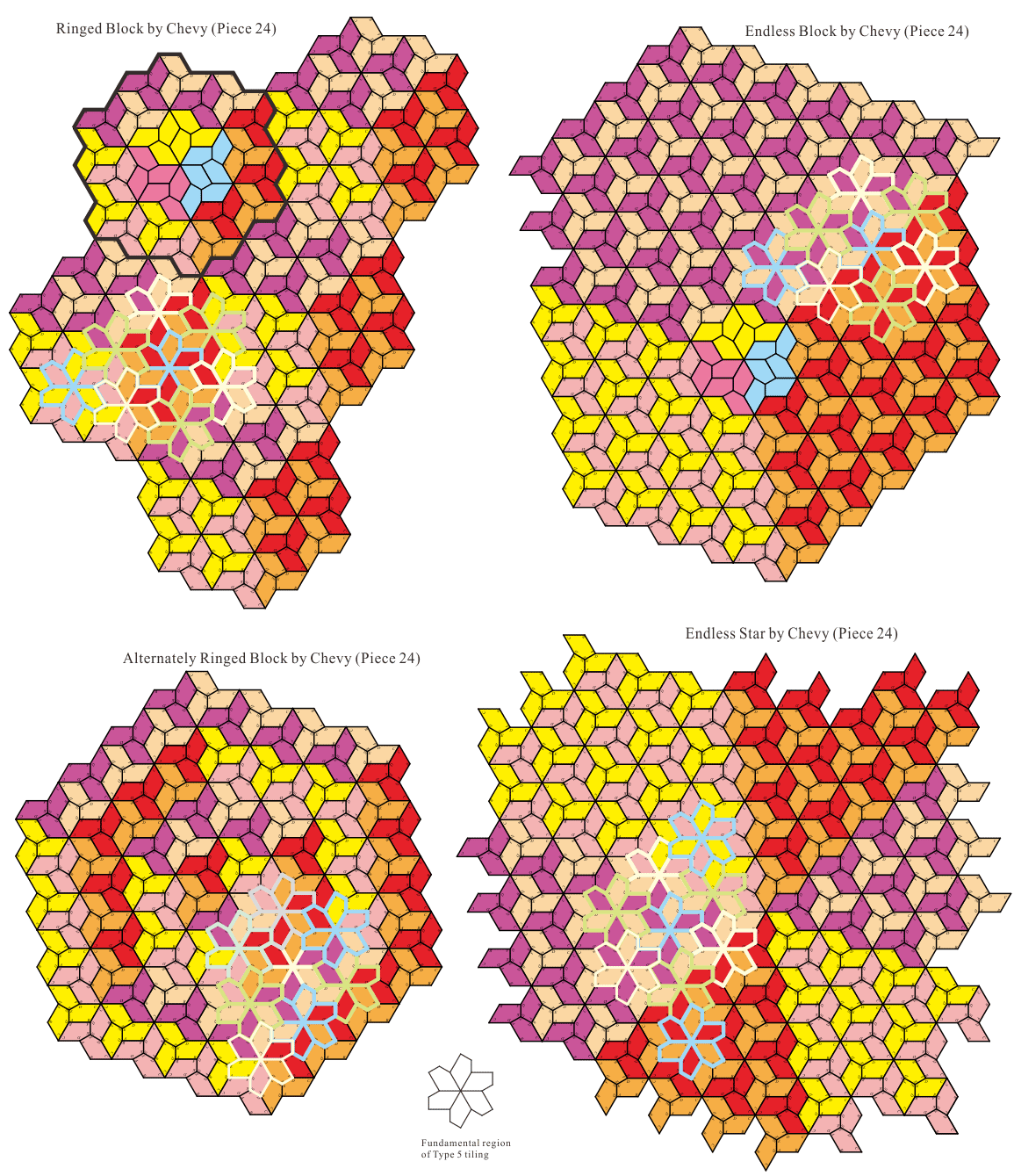} 
  \caption{{\small 
Examples of Block that converted to convex pentagon tiling.} 
\label{fig4}
}
\end{figure}

The tilings shown in Section 8.5 of 37sec824 and 37sec800 correspond to 
convex pentagon tiling of Type 5 (see Figure 2 in \cite{S_and_A_2017}), 
Rice1995-tiling (see Figure 4 in \cite{S_and_A_2017}), Figure 25 in 
\cite{S_and_A_2017}, etc. If the tilings of Zigzag and Block shown in Section 8.2 of 
37sec800 are filled with anterior TH-pentagons, then they are Type 5 tiling 
as shown in Figures~\ref{fig3} and \ref{fig4}. As mentioned in 
\cite{S_and_A_2017}, since the HFL1-unit can be reversed freely, it will 
also seem that they are the convex pentagon tilings combining 
AHFL1-units and PHFL1-units (see Figure 16 in \cite{S_and_A_2017}).

\begin{figure}[htbp]
 \centering\includegraphics[width=15cm,clip]{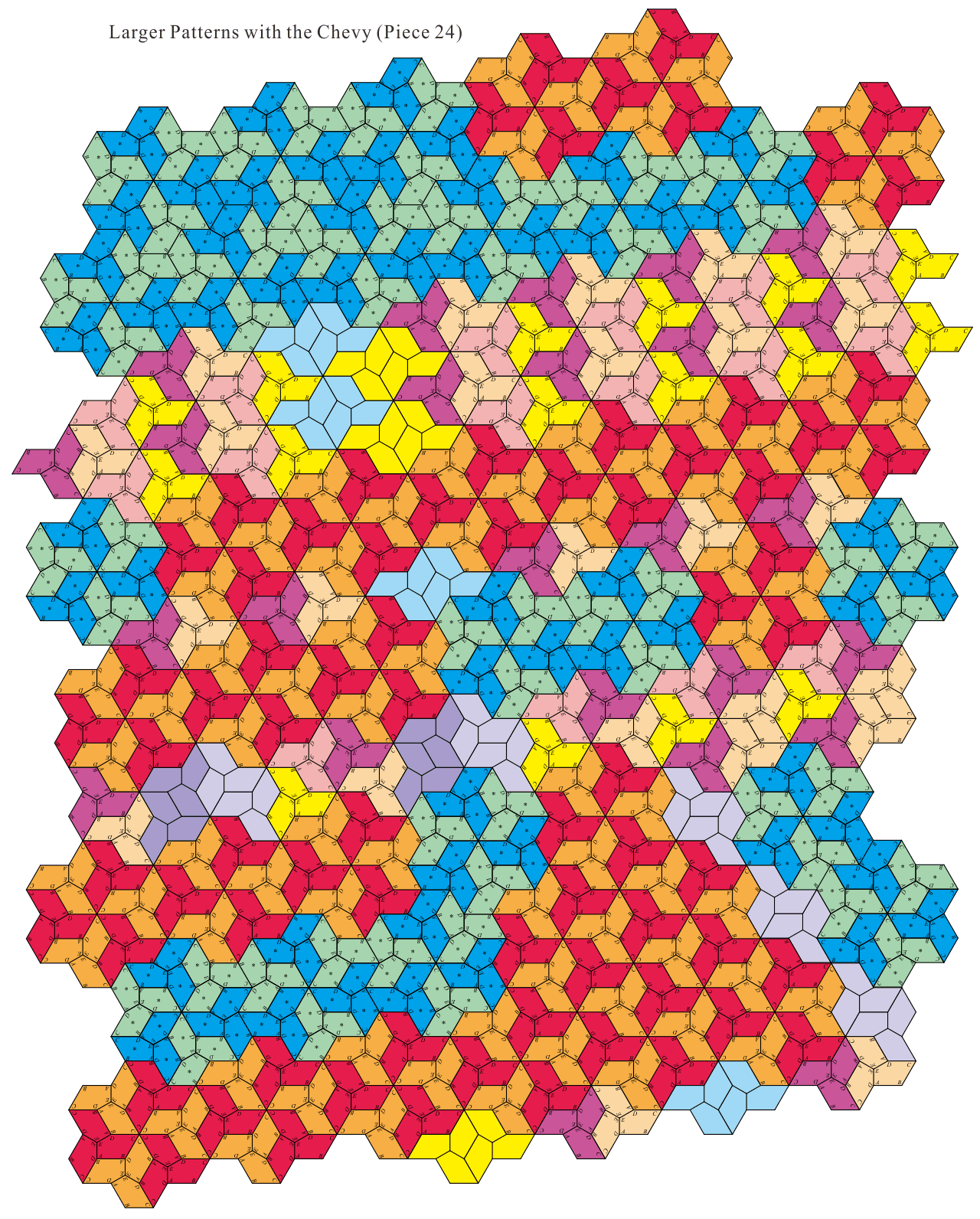} 
  \caption{{\small 
Example of Larger Patterns with Chevy that converted to convex 
pentagon tiling.} 
\label{fig5}
}
\end{figure}

\begin{figure}[htbp]
 \centering\includegraphics[width=15cm,clip]{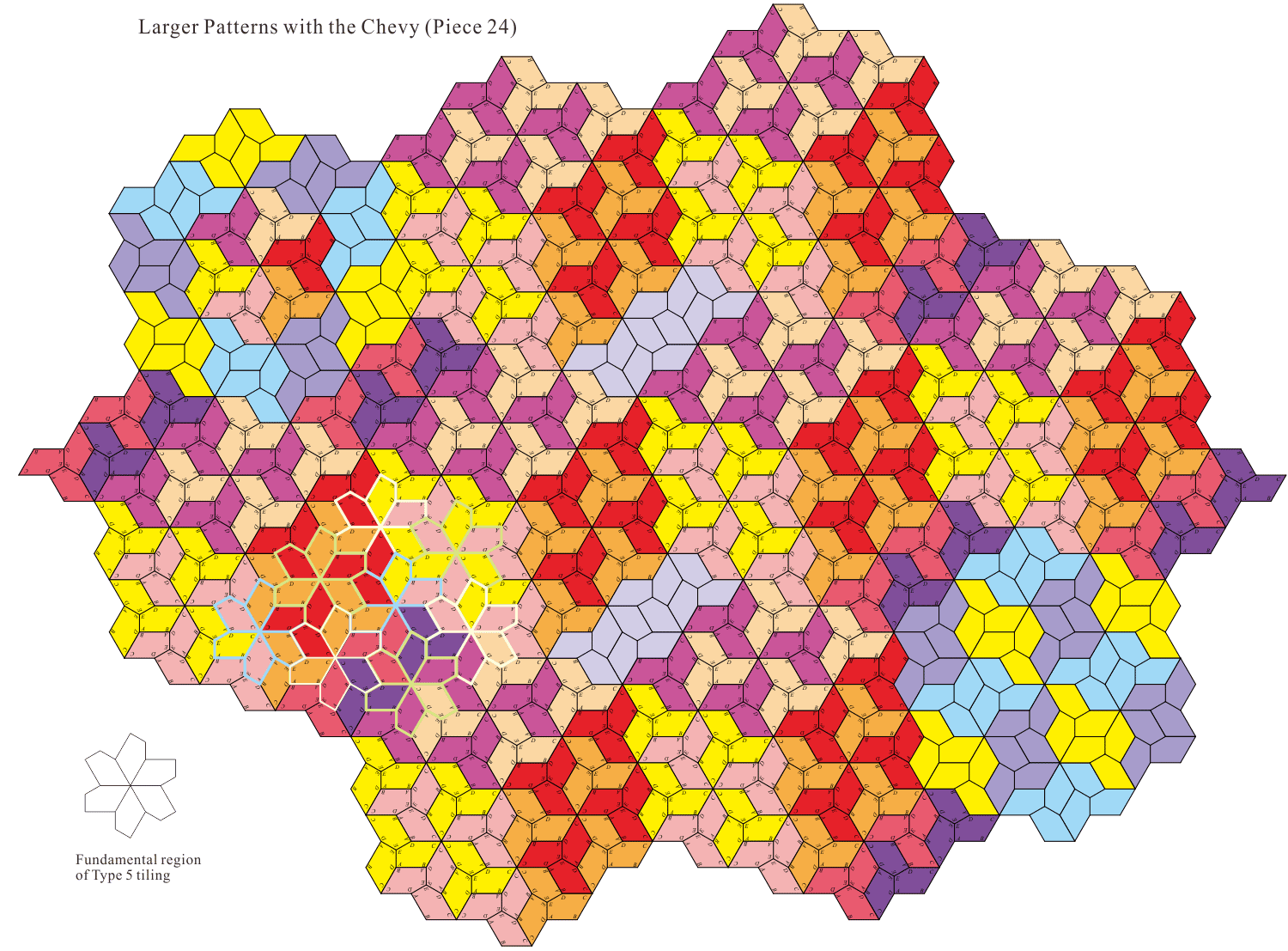} 
  \caption{{\small 
Example of Larger Patterns with Chevy that converted to convex 
pentagon tiling.} 
\label{fig6}
}
\end{figure}

\begin{figure}[htbp]
 \centering\includegraphics[width=15cm,clip]{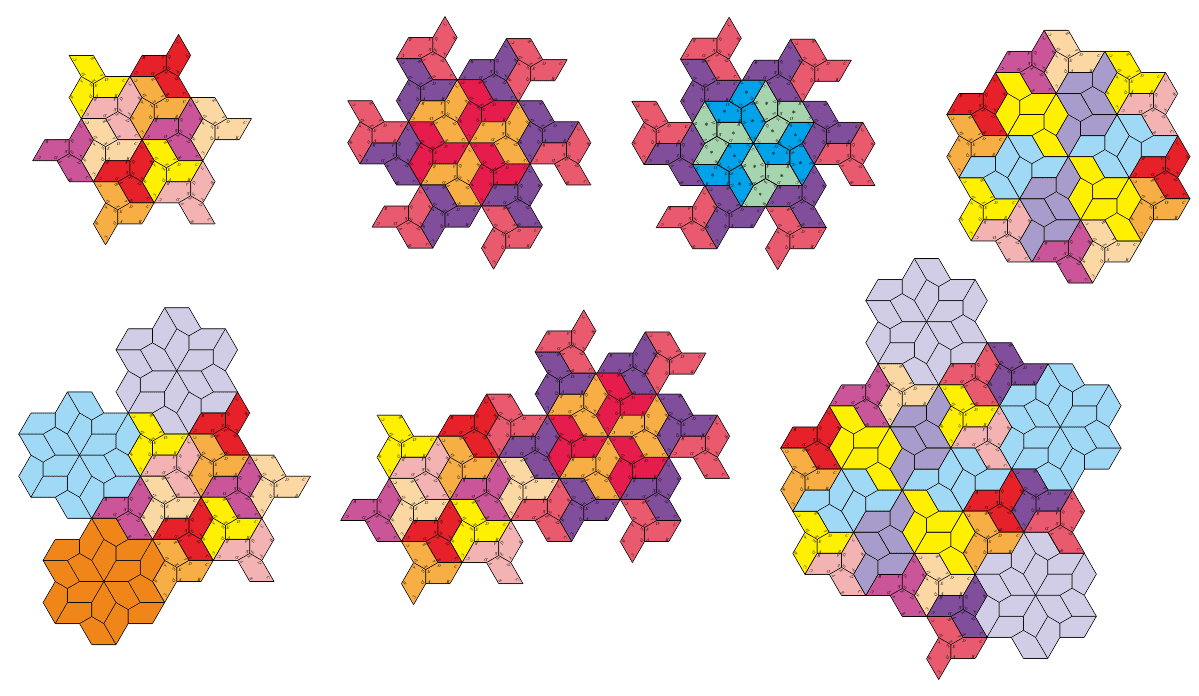} 
  \caption{{\small 
Examples of Larger Patterns with Chevy that converted to convex 
pentagon tiling.} 
\label{fig7}
}
\end{figure}

\begin{figure}[htbp]
 \centering\includegraphics[width=15cm,clip]{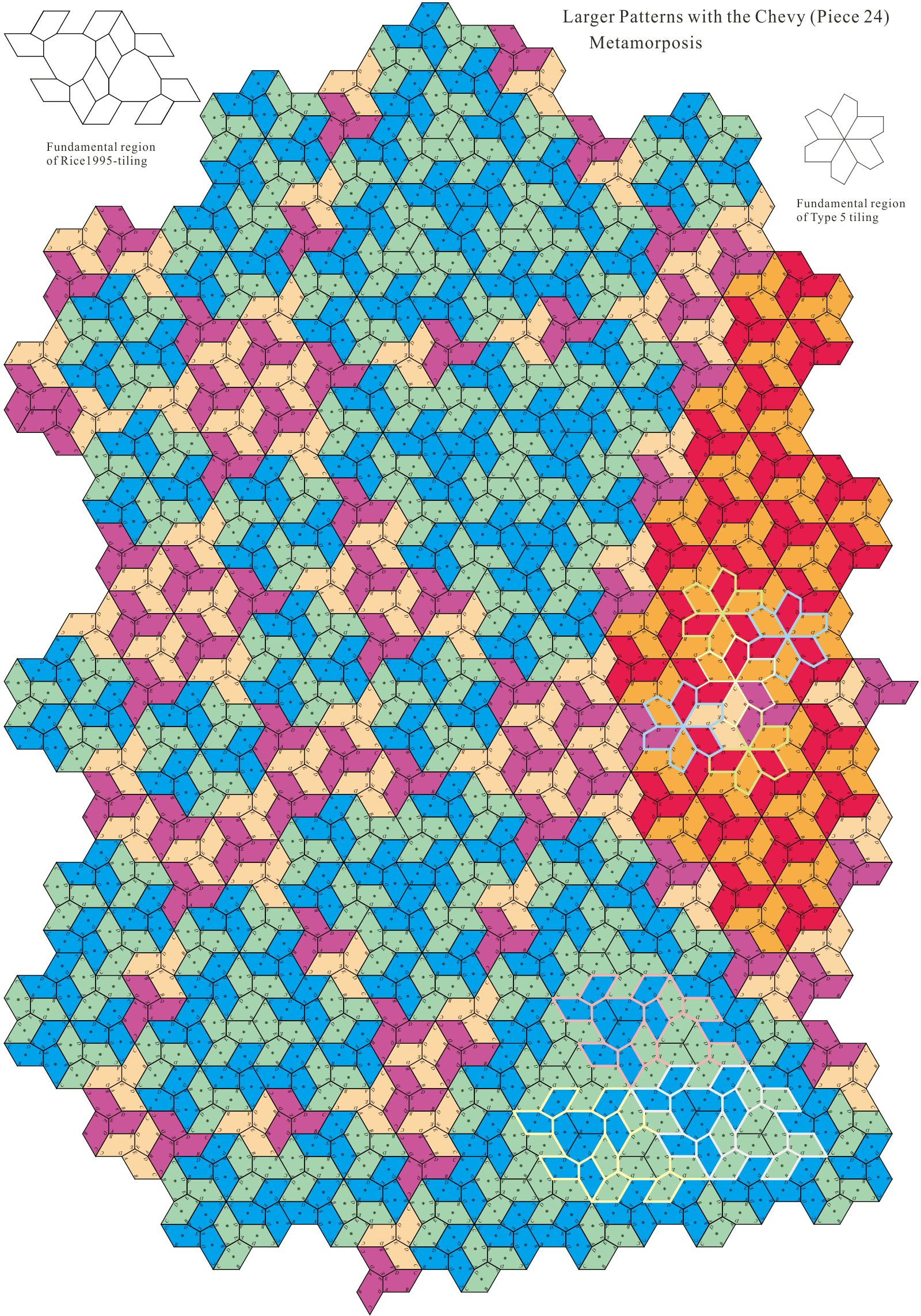} 
  \caption{{\small 
Example of Metamorposis that converted to convex pentagon tiling.} 
\label{fig8}
}
\end{figure}

The tilings shown in Larger Patterns in Section 8.2 of 37sec800 correspond 
to known convex pentagon tilings or the convex pentagon tilings combining 
anterior windmill units, posterior windmill units, AHFL1-units, and 
PHFL1-units, if Chevies and Miniflowers in the tilings are filled with 
TH-pentagon (see Figures~\ref{fig5}, \ref{fig6}, and \ref{fig7}). The tilings of 
each Stage that are called the last Metamorposis of Section 8.3 in 
37sec800 correspond to tilings that can be formed according to Case 2 
in Section 3 in \cite{S_and_A_2017} (see Figure ~\ref{fig8}). If Chevies and 
Miniflowers in the tilings of each Stage are 
filled with anterior TH-pentagon, they are Type 5 tilings.

\renewcommand{\figurename}{{\small Figure.}}

\renewcommand{\figurename}{{\small Figure.}}

\renewcommand{\figurename}{{\small Figure.}}

\renewcommand{\figurename}{{\small Figure.}}

\renewcommand{\figurename}{{\small Figure.}}

\renewcommand{\figurename}{{\small Figure.}}

\section{Tilings by only the ship units}
\label{section3}

In this section, the tiling of Hindriks' site \cite{Hindriks-site} are introduced 
which corresponds to tiling by only the ship unit. The articles related to 
Heptiamond 15 on Hindriks' site \cite{Hindriks-site} are as follows.

\medskip

37sec815 $\colon$ \url{http://www.jhhindriks.info/37/37sec815.htm}

37sec800 $\colon$ Section 8.4 in \url{http://www.jhhindriks.info/37/37sec800.htm}

\medskip

Hindriks calls the figure which corresponds to the hexagonal flower L2 unit 
(HFL2-unit) in \cite{S_and_A_2017} a Flower. The Flower by Heptiamond 
15 corresponds to Figure 41(n) in \cite{S_and_A_2017}.

Hereafter, it shows Figure of \cite{S_and_A_2017} when there are tilings 
corresponding to Hindriks' tilings already, otherwise it shows figures 
converted to convex pentagonal tilings. For Classes S1, S2, S3, S4, and 
S5, refer to Section 4 (Figure 31) in \cite{S_and_A_2017}.

Rotational Pair 2 of 37sec815 correspond to the tiling with Class S1 of 
Figure 32 in \cite{S_and_A_2017}.

\begin{figure}[htbp]
 \centering\includegraphics[width=15cm,clip]{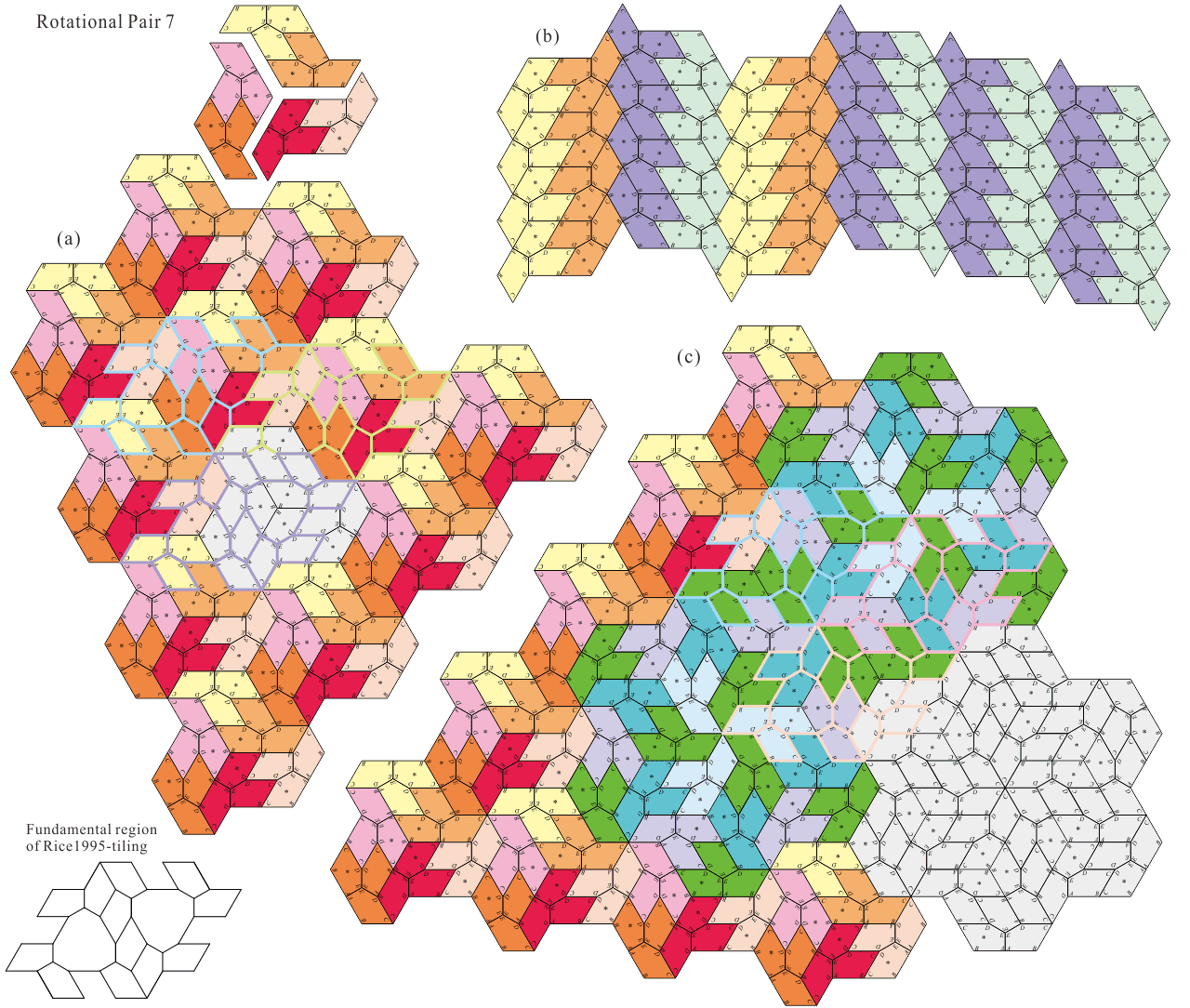} 
  \caption{{\small 
Tilings with Class S2.} 
\label{fig9}
}
\end{figure}

Rotational Pair 7 of 37sec815 correspond to the tiling with Class S2. 
Figure 33 in \cite{S_and_A_2017} is contained in this pair. As shown in 
Figure~\ref{fig9}(a), the tiling using the form of three Rotation Pairs 7 
indicated by Hindriks as the fundamental region is the Rice1995-tiling 
formed only by the ship unit (see Figure 56(b) in \cite{S_and_A_2017}). 
On the other hand, tiling (combination) as shown in Figure~\ref{fig9}(b) 
is a case unknown to the authors. Tilings including the Flower of 
Rotational Pair 7 as shown in Figure~\ref{fig9}(c) correspond to the 
tilings by HFL2-units of Figure 41(n) in \cite{S_and_A_2017} (see 
Figure 56 in \cite{S_and_A_2017}), the Rice1995-tiling by only the ship 
units, or cases that connected those two tilings.

Rotational Pair 9 of 37sec815 correspond to the tiling with Class S3 of 
Figure 34 in \cite{S_and_A_2017} or the tiling with Classes S1 and S3 of 
Figure 37 in \cite{S_and_A_2017}.

Crooked Pair of 37sec815 correspond to the tiling with Class S4 in 
Figures 35 and 36 in \cite{S_and_A_2017}.

\begin{figure}[htbp]
 \centering\includegraphics[width=15cm,clip]{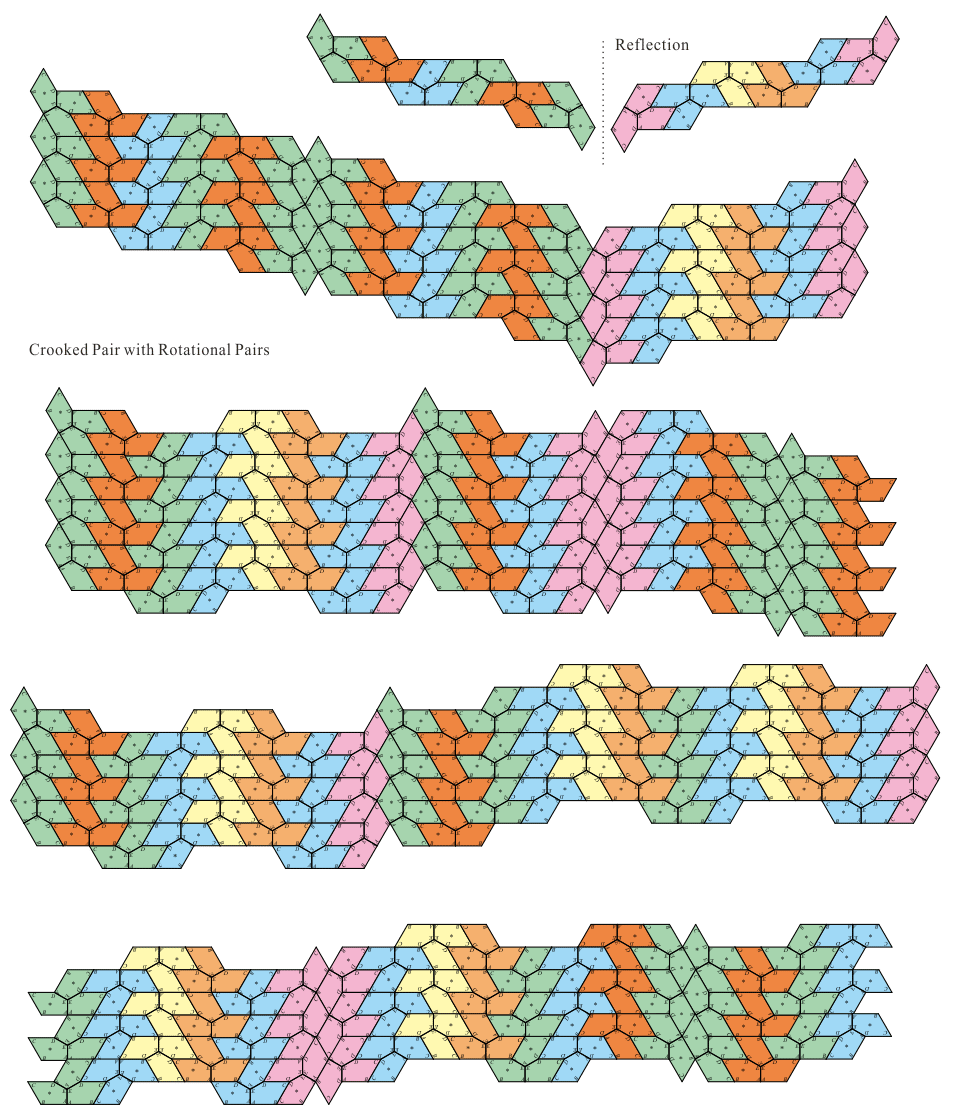} 
  \caption{{\small 
Tilings with Classes S1 and S4.} 
\label{fig10}
}
\end{figure}

\begin{figure}[htbp]
 \centering\includegraphics[width=14.5cm,clip]{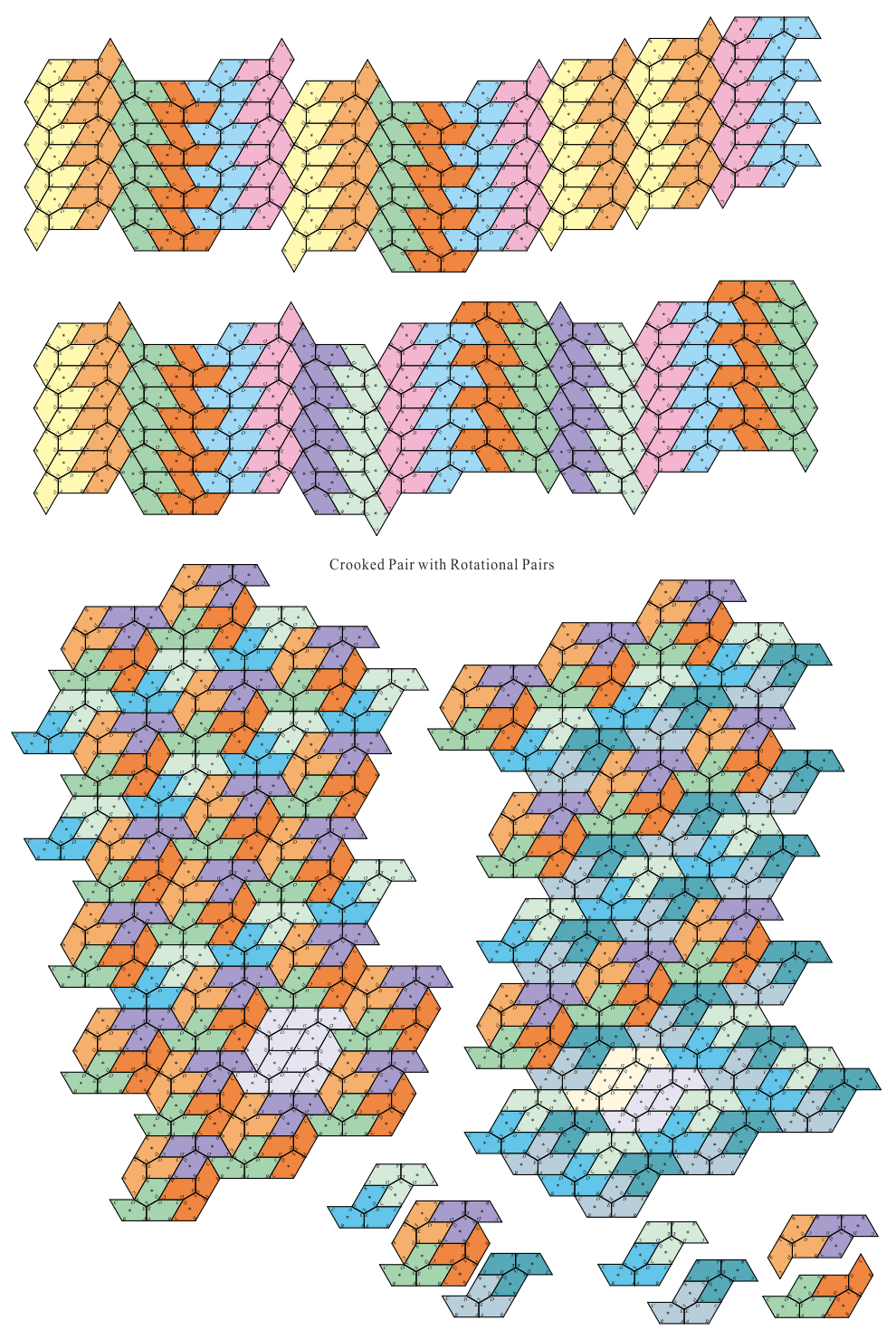} 
  \caption{{\small 
Tilings with Classes S2 and S4.} 
\label{fig11}
}
\end{figure}

\begin{figure}[htbp]
 \centering\includegraphics[width=15cm,clip]{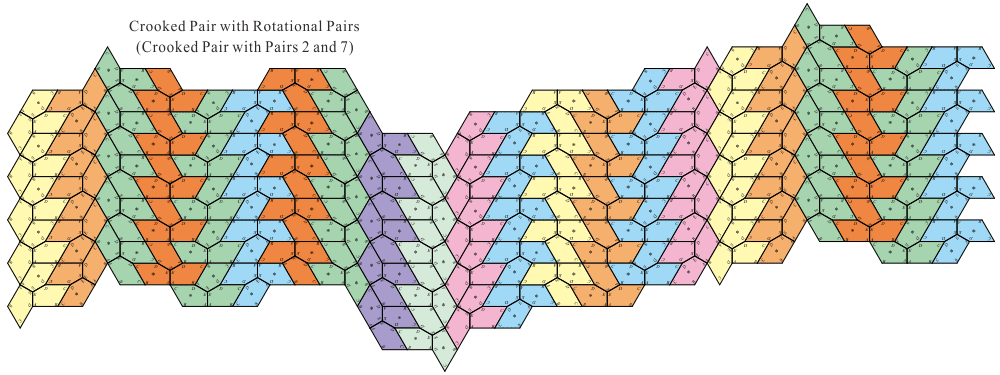} 
  \caption{{\small 
Tiling with Classes S1, S2, and S4.} 
\label{fig12}
}
\end{figure}

Crooked Pair with Rotational Pairs of 37sec815 contain the tiling with 
Classes S2 and S4 of Figure 38 in \cite{S_and_A_2017}. For the tilings 
corresponding to this case, the authors have found only the tilings of 
Figure 38 in \cite{S_and_A_2017}. For other tilings not known to the 
authors, it is confirmed that there is a tiling with Classes S1 and S4 
(see Figure~\ref{fig10}), a tiling with Classes S2 and S4 (see 
Figure~\ref{fig11}), and a tiling with Classes S1, S2 and S4 (see Figure~\ref{fig12}).

\begin{figure}[htbp]
 \centering\includegraphics[width=15cm,clip]{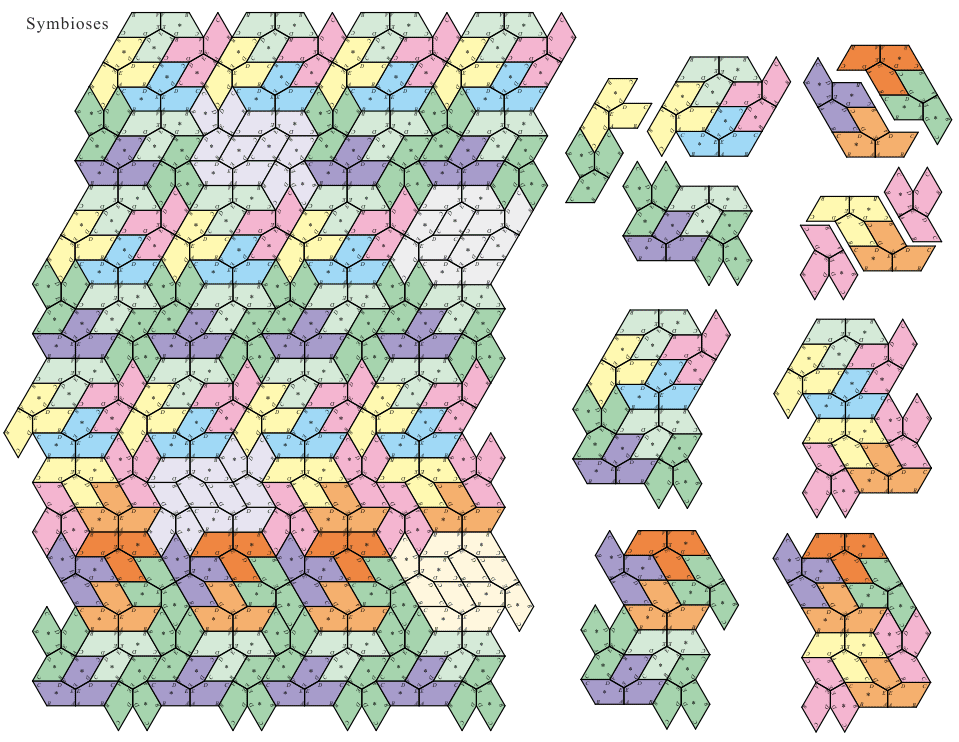} 
  \caption{{\small 
Tiling with Classes S2, S4, and S5.} 
\label{fig13}
}
\end{figure}

Symbiose of 37sec815 are tilings with Classes S2, S4, and S5 
(see Figure~\ref{fig13}). They are tilings unknown to the authors.

\begin{figure}[htbp]
 \centering\includegraphics[width=15cm,clip]{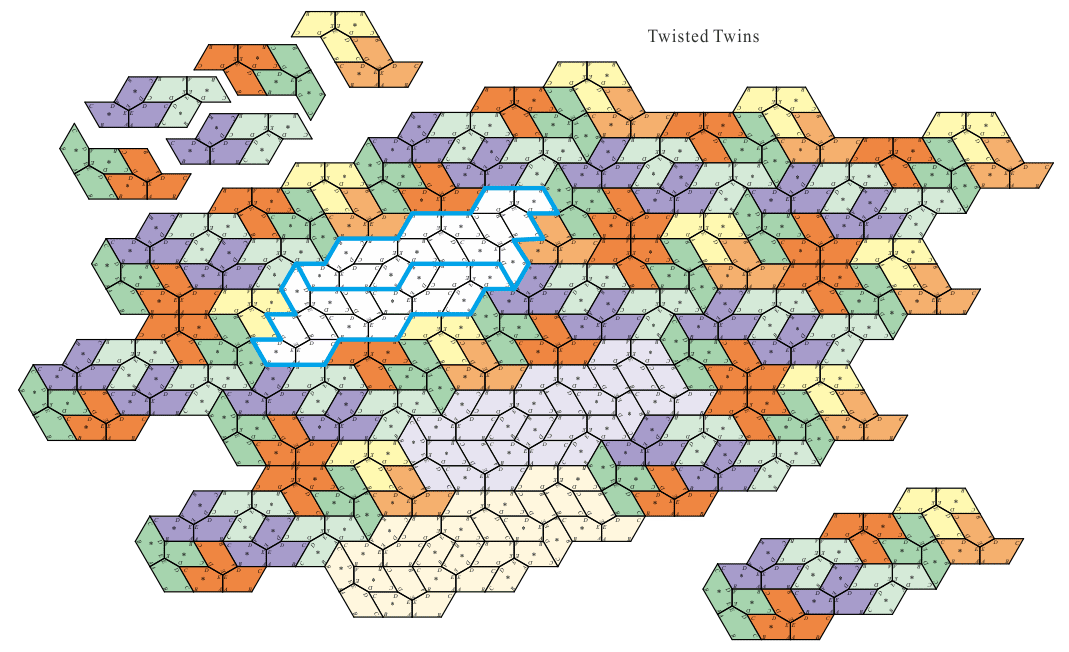} 
  \caption{{\small 
Tiling with Classes S1, S2, and S4.} 
\label{fig14}
}
\end{figure}

\begin{figure}[htbp]
 \centering\includegraphics[width=15cm,clip]{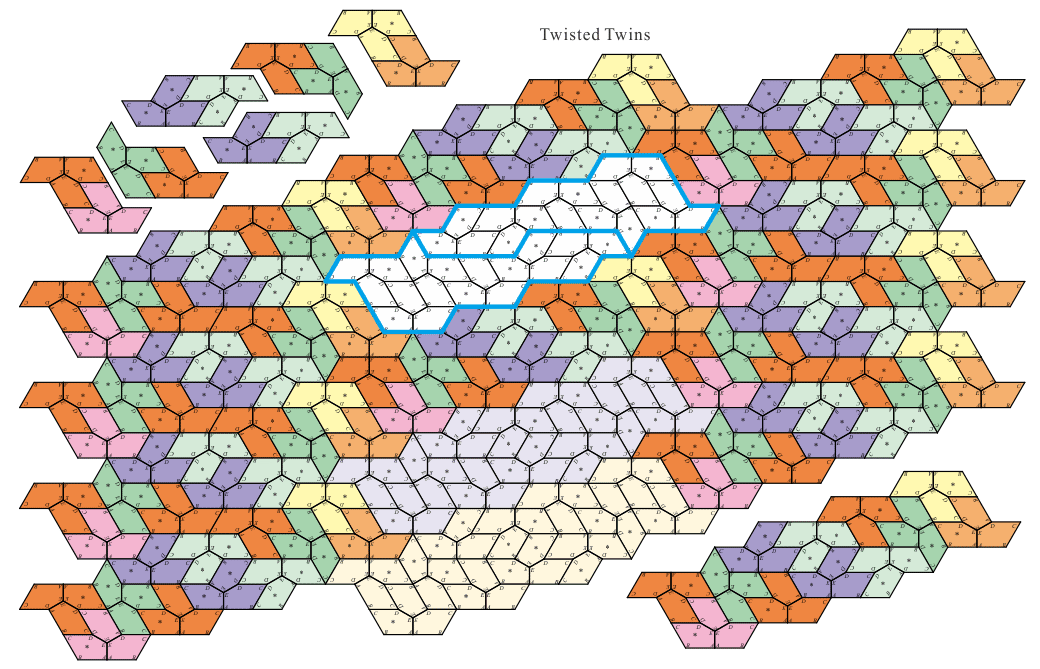} 
  \caption{{\small 
Tiling with Classes S1, S2, and S4.} 
\label{fig15}
}
\end{figure}

\begin{figure}[htbp]
 \centering\includegraphics[width=15cm,clip]{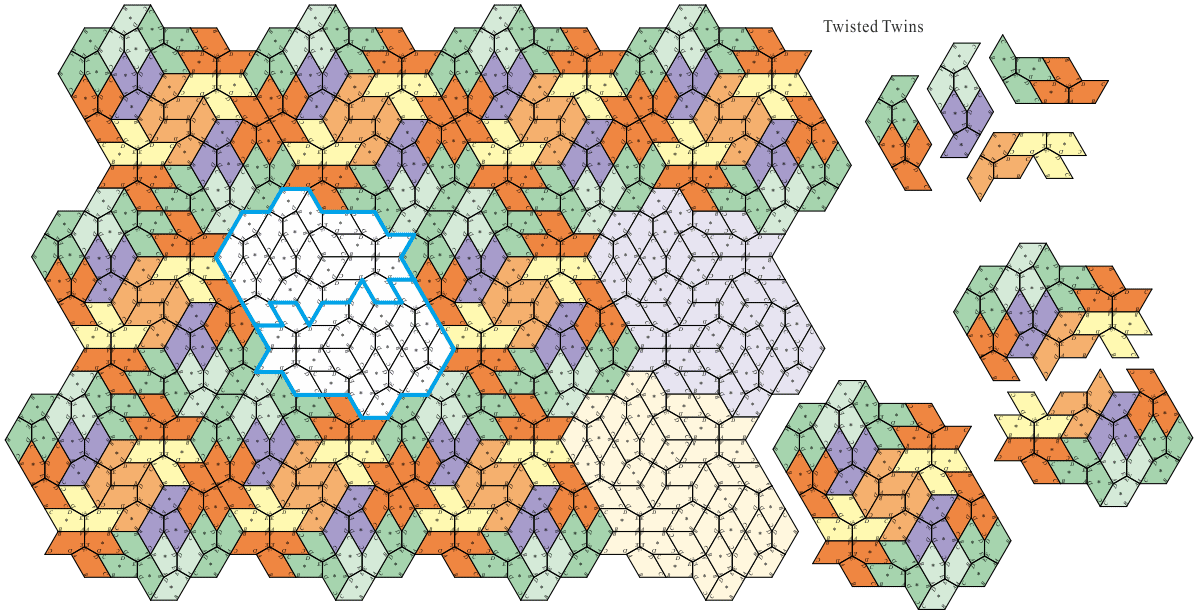} 
  \caption{{\small 
Tiling with Classes S2, S3, and S4.} 
\label{fig16}
}
\end{figure}

\begin{figure}[htbp]
 \centering\includegraphics[width=15cm,clip]{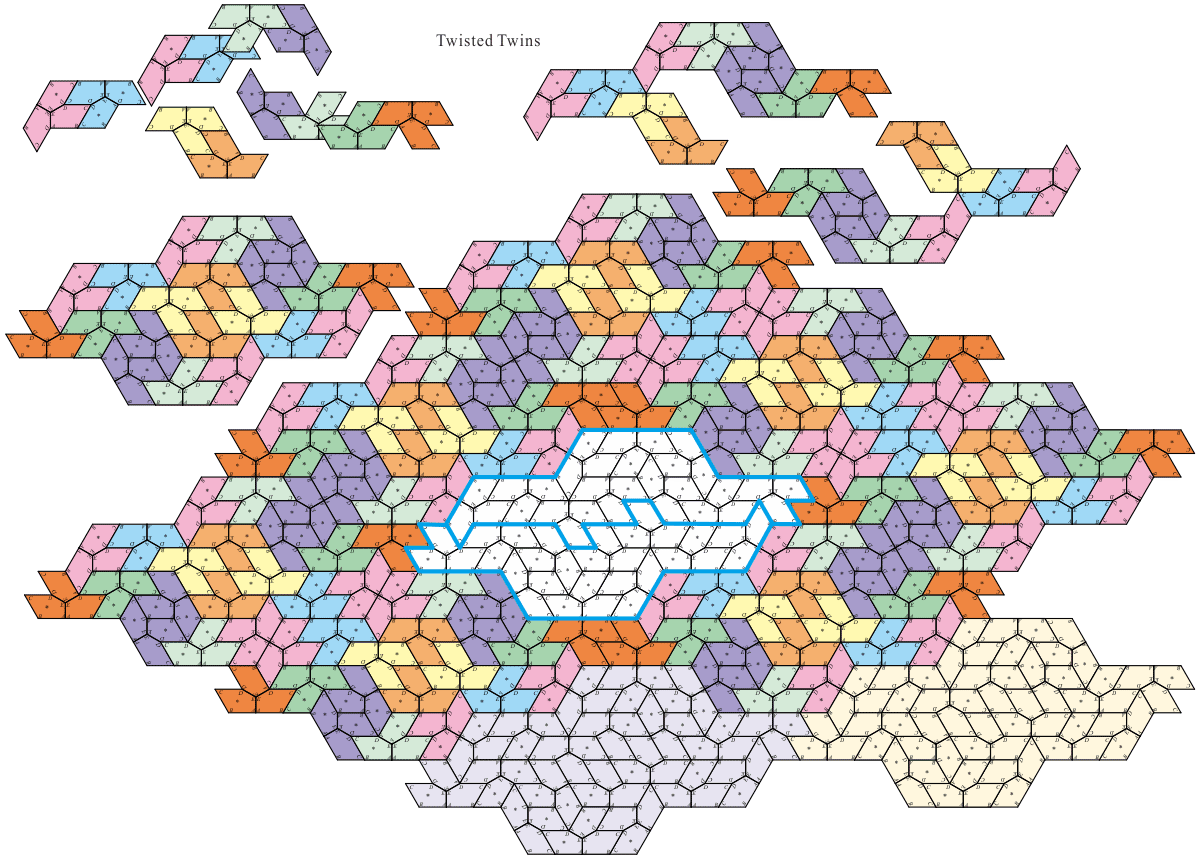} 
  \caption{{\small 
Tiling with Classes S2, S3, S4, and S5.} 
\label{fig17}
}
\end{figure}

\begin{figure}[htbp]
 \centering\includegraphics[width=15cm,clip]{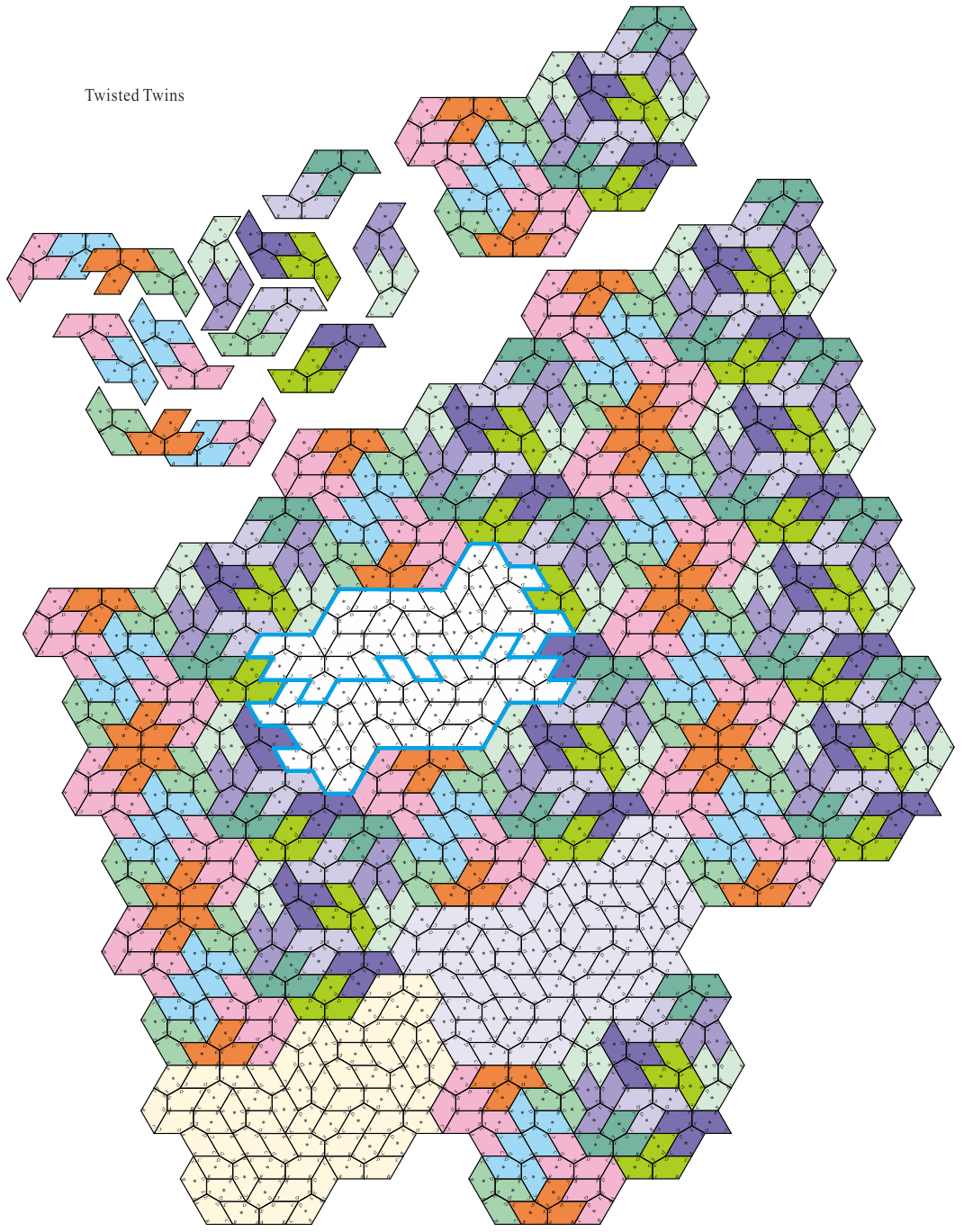} 
  \caption{{\small 
Tiling with Classes S2, S3, and S4.} 
\label{fig18}
}
\end{figure}

\begin{figure}[htbp]
 \centering\includegraphics[width=15cm,clip]{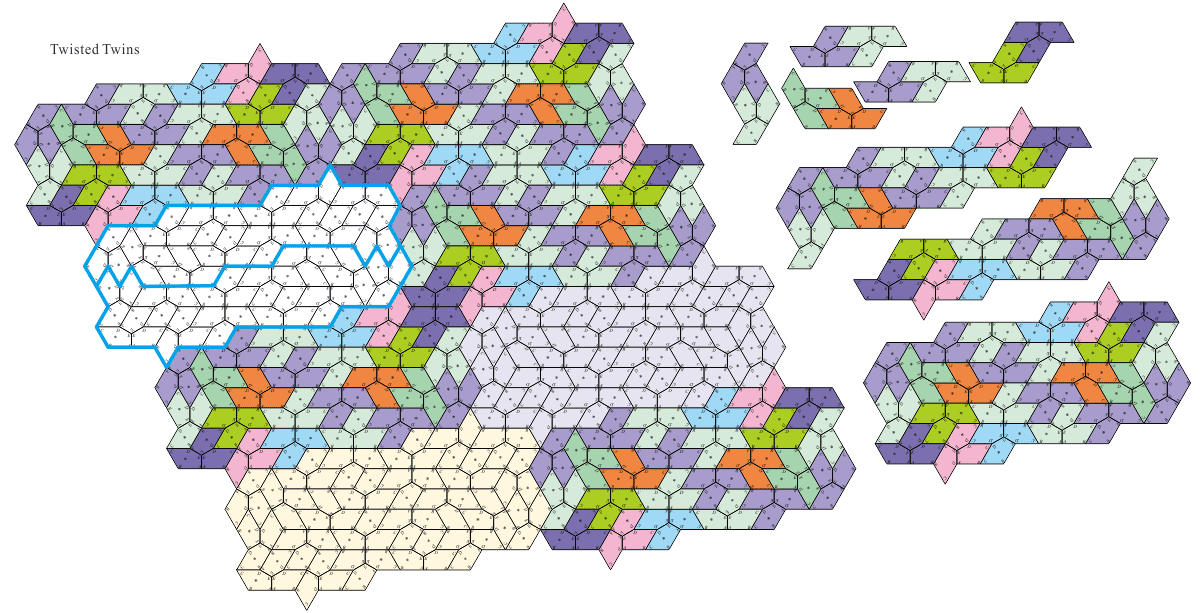} 
  \caption{{\small 
Tiling with Classes S1, S2, and S4.} 
\label{fig19}
}
\end{figure}

\begin{figure}[htbp]
 \centering\includegraphics[width=15cm,clip]{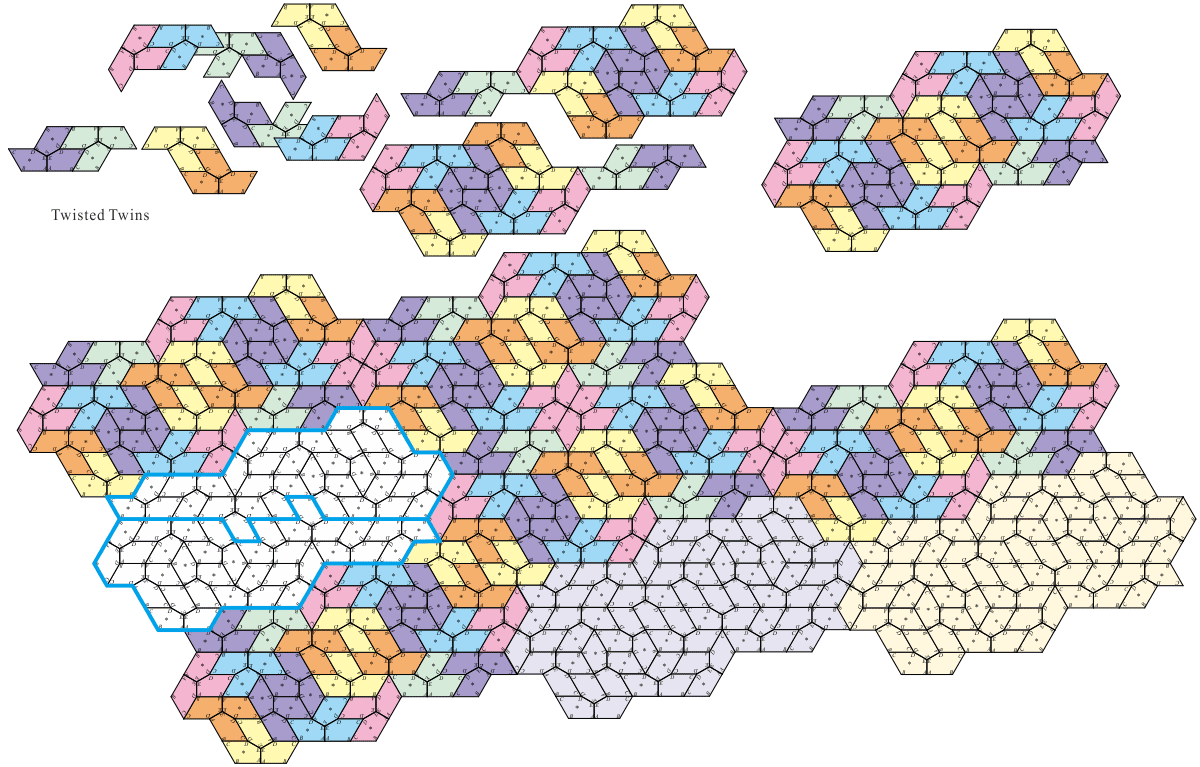} 
  \caption{{\small 
Tiling with Classes S1, S2, S3, and S4.} 
\label{fig20}
}
\end{figure}

\begin{figure}[htbp]
 \centering\includegraphics[width=15cm,clip]{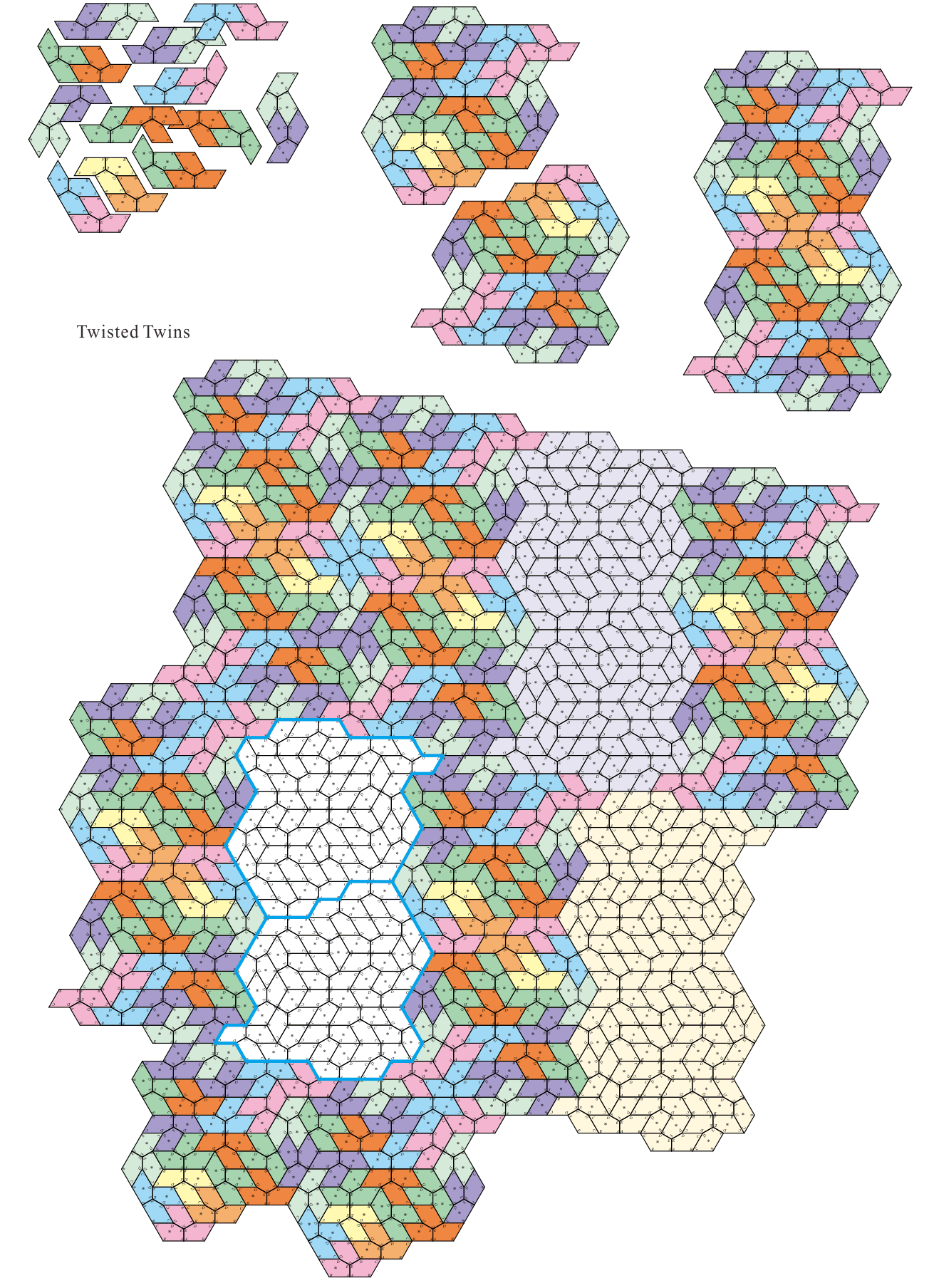} 
  \caption{{\small 
Tiling with Classes S1, S2, S3, S4, and S5.} 
\label{fig21}
}
\end{figure}

\begin{figure}[htbp]
 \centering\includegraphics[width=15cm,clip]{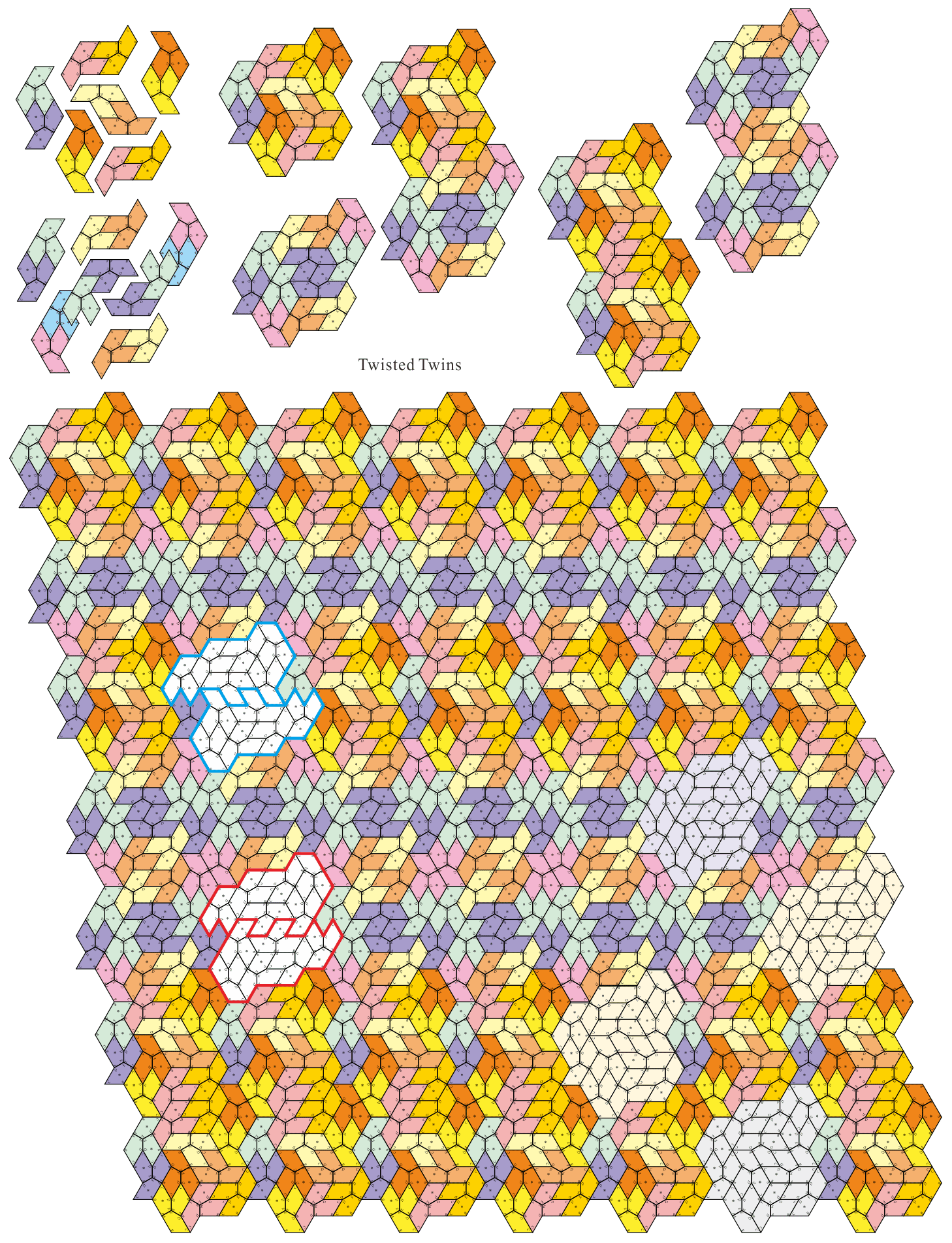} 
  \caption{{\small 
Tiling with Classes S1, S2, S3, and S4.} 
\label{fig22}
}
\end{figure}

Twisted Twins of 37sec815 are all tilings unknown to the authors (see 
Figures~\ref{fig14}-\ref{fig22}). For the classes used in each tiling, see 
the caption of each figure.

\begin{figure}[htbp]
 \centering\includegraphics[width=15cm,clip]{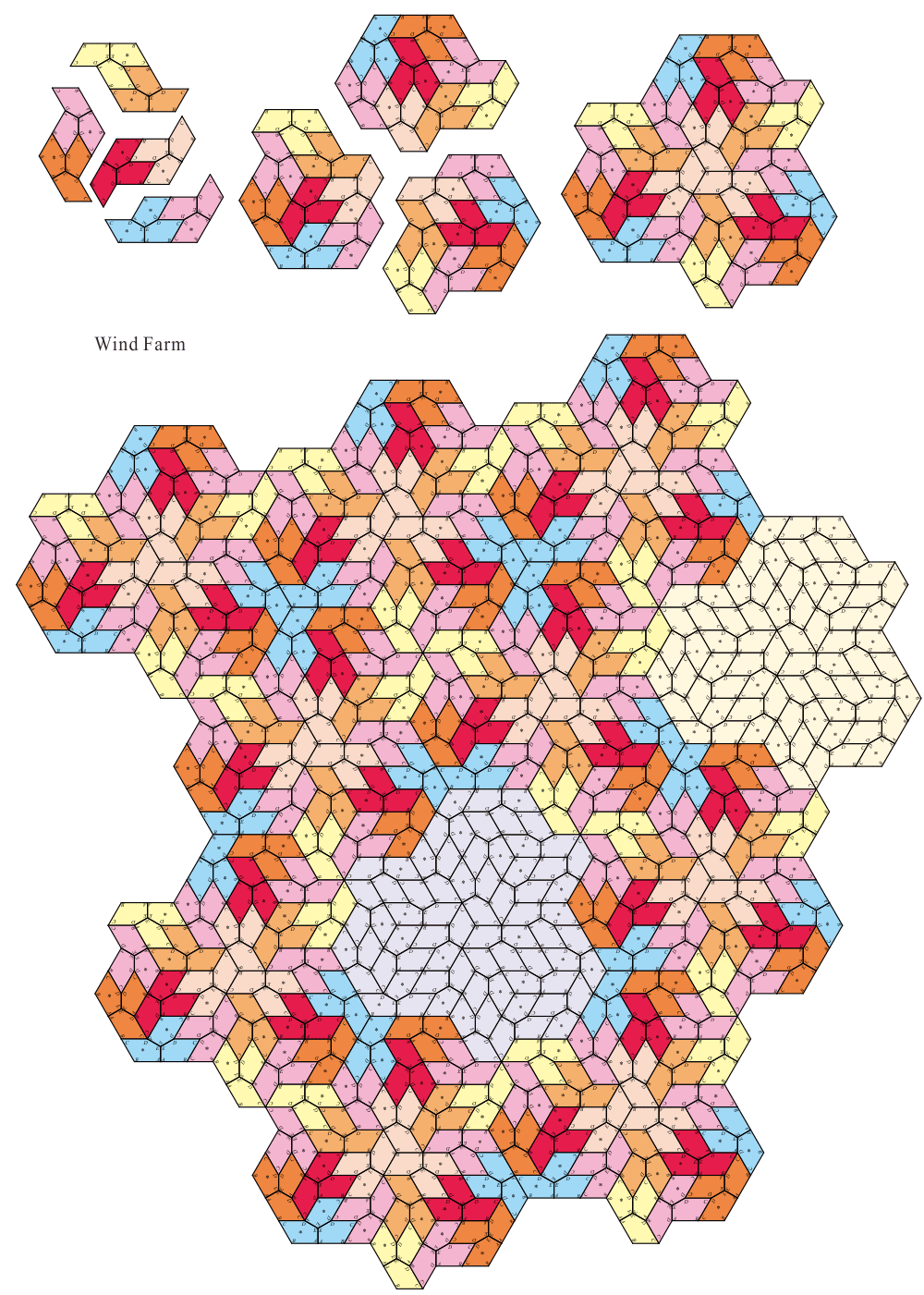} 
  \caption{{\small 
Tiling with Classes S2 and S4.} 
\label{fig23}
}
\end{figure}

 Wind Farm of 37sec815 is a tiling with Classes S2 and S4 (see Figure~\ref{fig23}). 
This tiling corresponds to a tiling formed by HFL2-units and a property of 
CN-units (see Figure 5 in \cite{S_and_A_2017} or Figure~\ref{fig30}). In other words, 
it is one variation of tilings that can be formed by the proposed method in 
\cite{S_and_A_2017}. Details are explained in the next section.

\begin{figure}[htbp]
 \centering\includegraphics[width=15cm,clip]{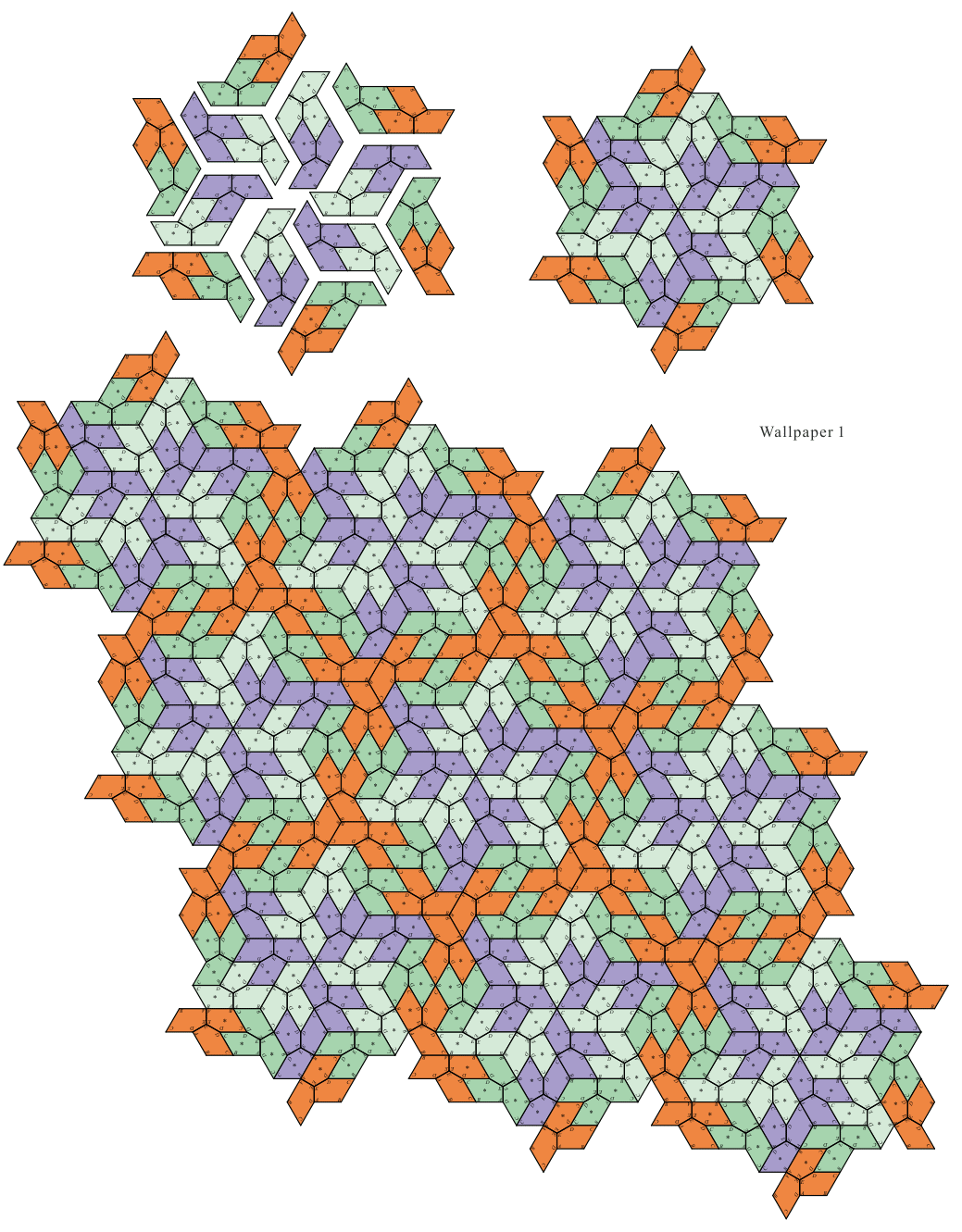} 
  \caption{{\small 
Tiling with Classes S2 and S4.} 
\label{fig24}
}
\end{figure}

 Wallpaper 1 of 37sec815 is a tiling with Classes S2 and S4 (see Figure~\ref{fig24}). 
This tiling corresponds to a tiling formed by HFL2-units and property of 
CN-units (see Figure 5 in \cite{S_and_A_2017}). Details are explained in the next section.

\begin{figure}[htbp]
 \centering\includegraphics[width=15cm,clip]{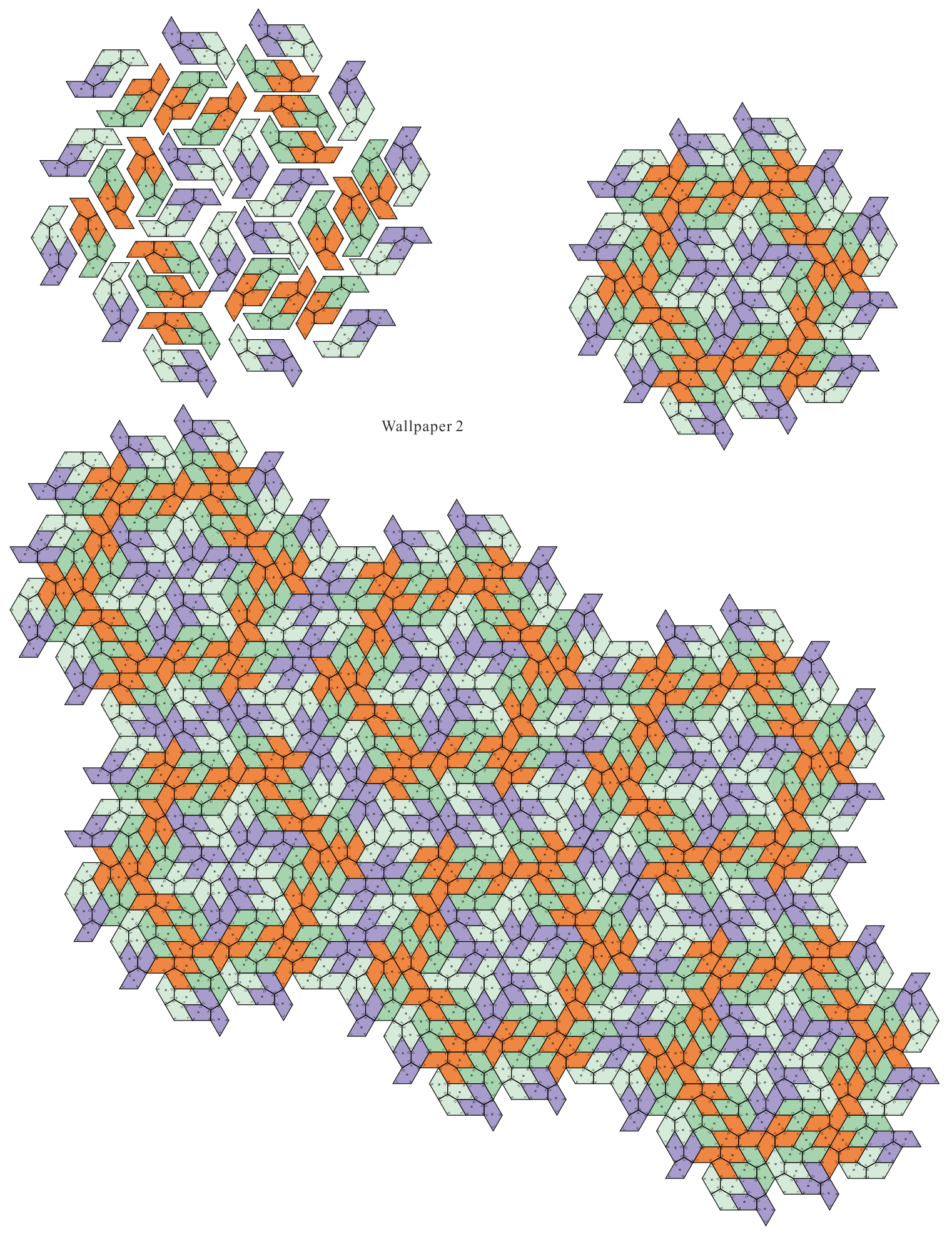} 
  \caption{{\small 
Tiling with Classes S2 and S4.} 
\label{fig25}
}
\end{figure}

Wallpaper 2 of 37sec815 is a tiling with Classes S2 and S4 (see Figure~\ref{fig25}). 
This tiling is a case unknown to the authors.

\begin{figure}[htbp]
 \centering\includegraphics[width=15cm,clip]{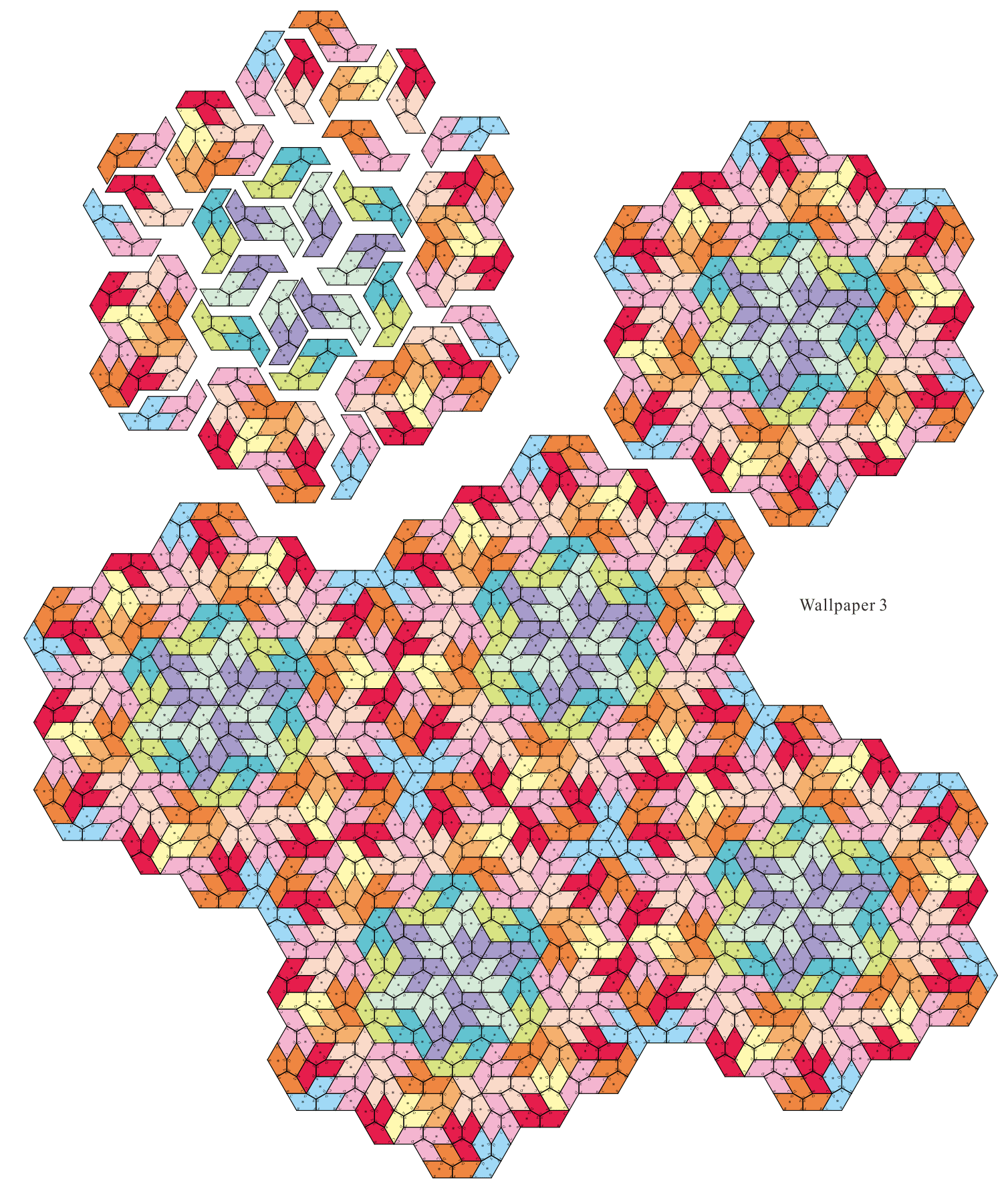} 
  \caption{{\small 
Tiling with Classes S2 and S4.} 
\label{fig26}
}
\end{figure}

Wallpaper 3 of 37sec815 is a tiling with Classes S2 and S4 (see Figure~\ref{fig26}). 
This tiling corresponds to a tiling formed by HFL2-units and a property of 
CN-units (see Figure 5 in \cite{S_and_A_2017}). Details are explained in the next section.

\begin{figure}[htbp]
 \centering\includegraphics[width=15cm,clip]{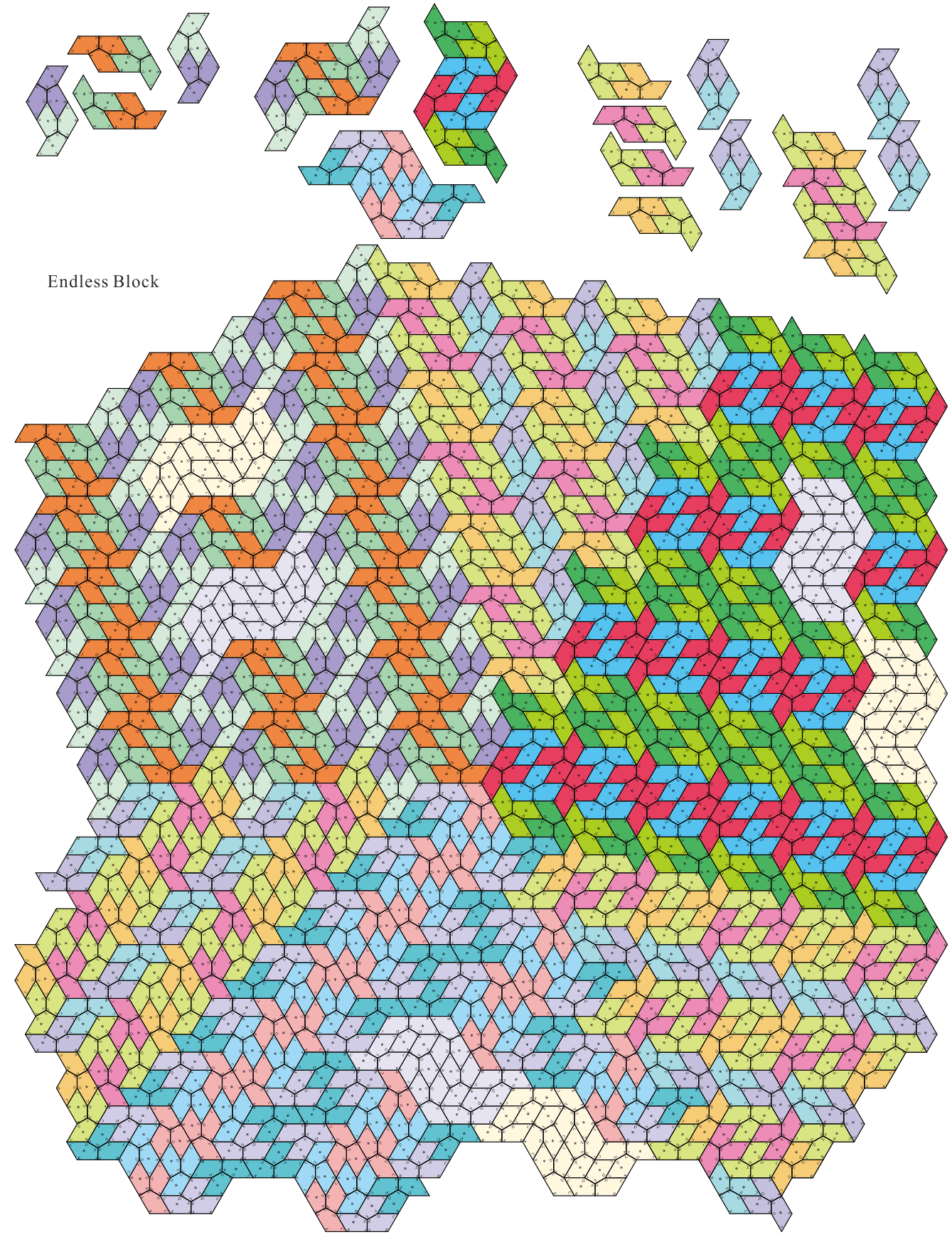} 
  \caption{{\small 
Tiling with Classes S2 and S4.} 
\label{fig27}
}
\end{figure}

Endless Block of 37sec815 is a tiling with Classes S2, S4, and S5 (see 
Figure~\ref{fig27}). This is formed by applying the property of tiling using Classes 
S2 and S4 in Figure 38 of \cite{S_and_A_2017}. This tiling is a case unknown to 
the authors.

\begin{figure}[htbp]
 \centering\includegraphics[width=15cm,clip]{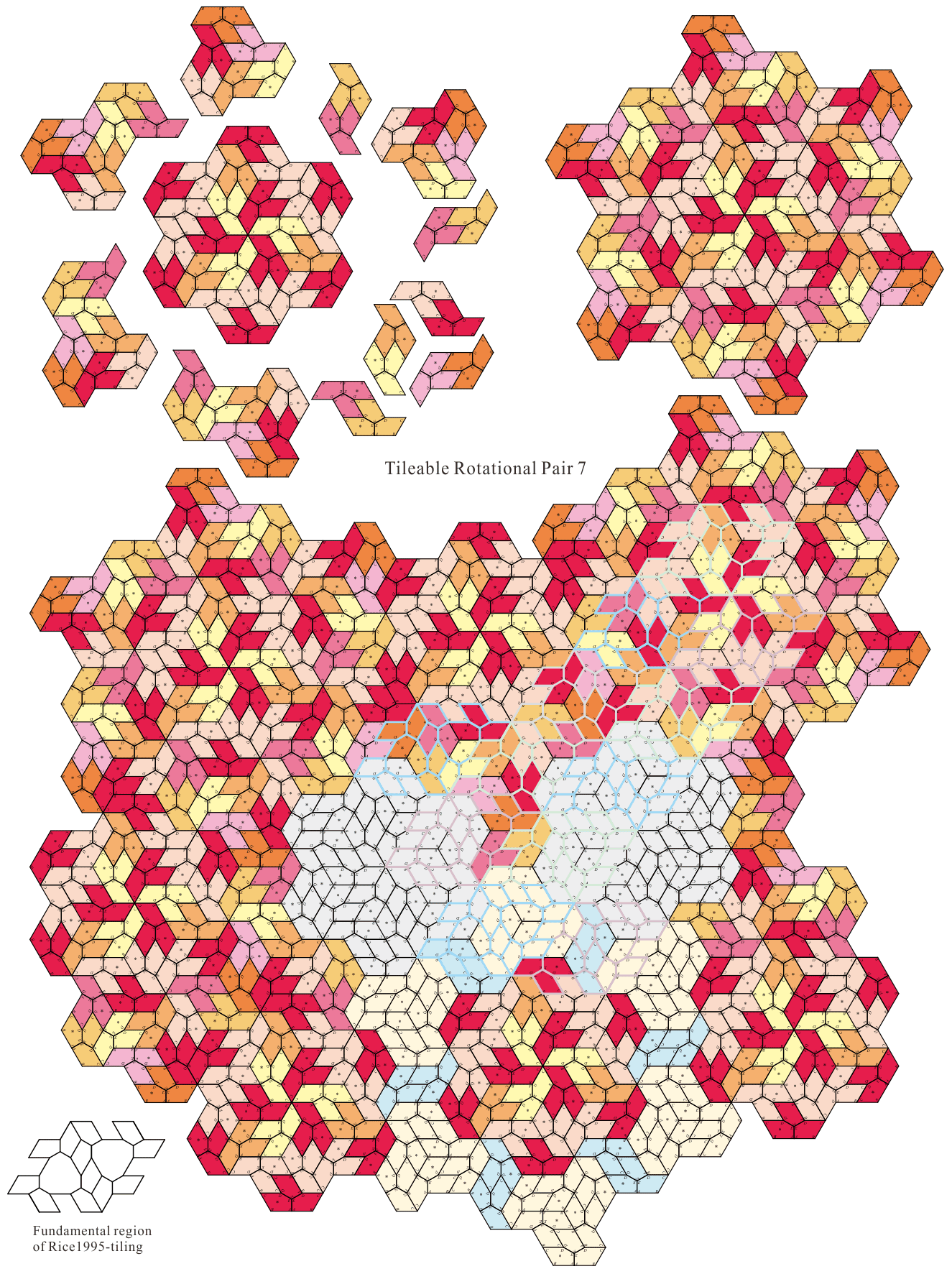} 
  \caption{{\small 
Tiling with Class S2.} 
\label{fig28}
}
\end{figure}

\begin{figure}[htbp]
 \centering\includegraphics[width=15cm,clip]{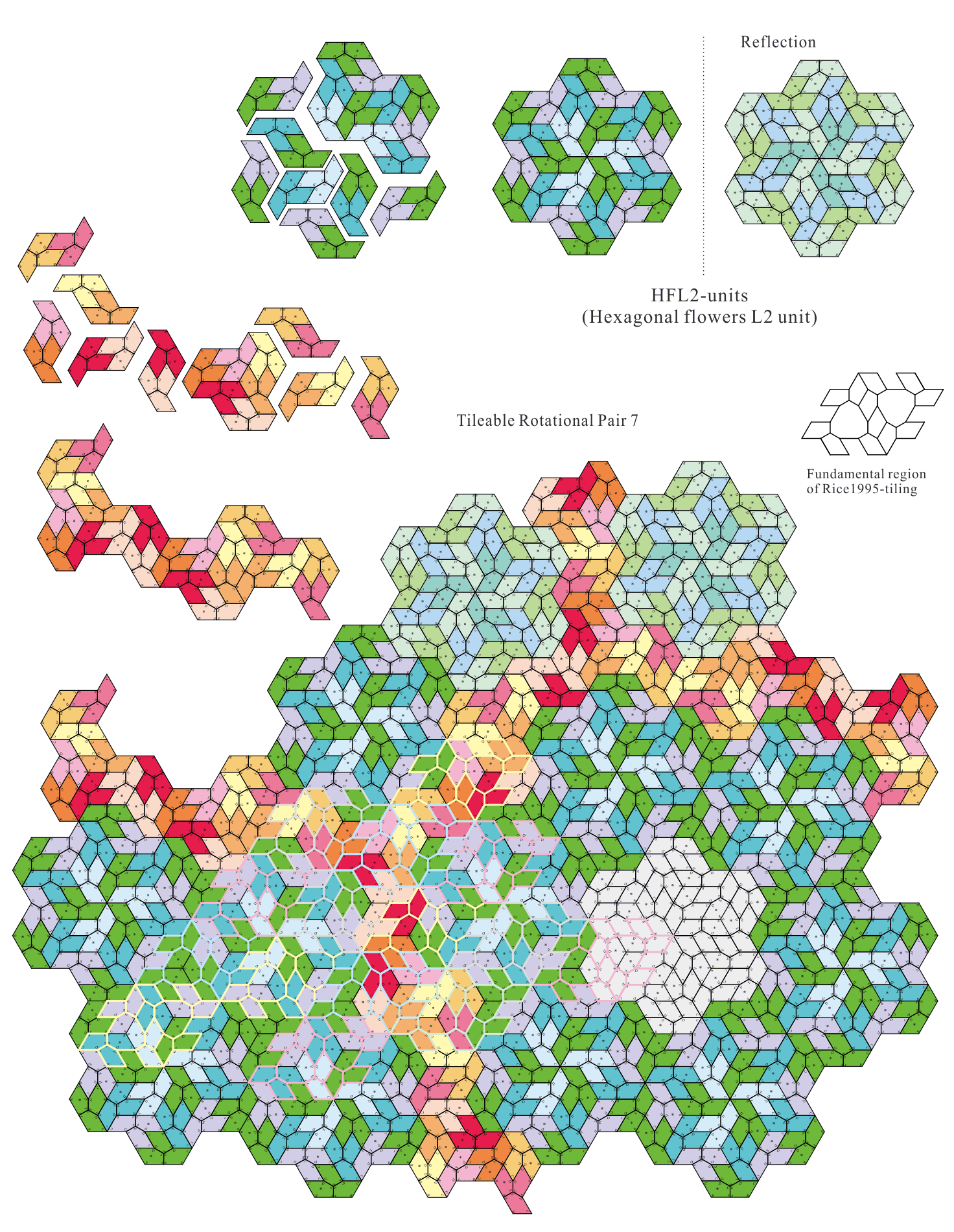} 
  \caption{{\small 
Tiling with Class S2.} 
\label{fig29}
}
\end{figure}

Tileable Rotational Pair 7 of Section 8.4 in 37sec800 are tilings with Class S2 
(see Figures~\ref{fig28} and \ref{fig29}). 
They are the same tilings as Figure~\ref{fig9}(c). 

\bigskip
\bigskip
\bigskip
\bigskip

\renewcommand{\figurename}{{\small Figure.}}

\bigskip
\bigskip
\bigskip
\bigskip
\bigskip
\bigskip
\bigskip

\renewcommand{\figurename}{{\small Figure.}}

\renewcommand{\figurename}{{\small Figure.}}

\renewcommand{\figurename}{{\small Figure.}}

\renewcommand{\figurename}{{\small Figure.}}

\renewcommand{\figurename}{{\small Figure.}}

\renewcommand{\figurename}{{\small Figure.}}

\renewcommand{\figurename}{{\small Figure.}}

\renewcommand{\figurename}{{\small Figure.}}

\renewcommand{\figurename}{{\small Figure.}}

\renewcommand{\figurename}{{\small Figure.}}

\renewcommand{\figurename}{{\small Figure.}}

\renewcommand{\figurename}{{\small Figure.}}

\renewcommand{\figurename}{{\small Figure.}}

\renewcommand{\figurename}{{\small Figure.}}

\renewcommand{\figurename}{{\small Figure.}}

\renewcommand{\figurename}{{\small Figure.}}

\renewcommand{\figurename}{{\small Figure.}}

\renewcommand{\figurename}{{\small Figure.}}

\renewcommand{\figurename}{{\small Figure.}}

\renewcommand{\figurename}{{\small Figure.}}

\section{Tilings whose convex nonagon units are reversed}
\label{section4}

\begin{figure}[htbp]
 \centering\includegraphics[width=10cm,clip]{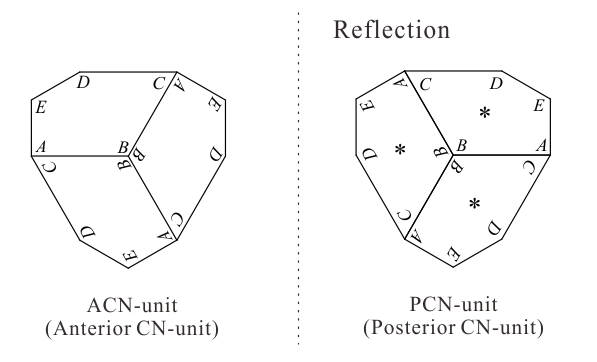} 
  \caption{{\small 
Convex nonagon unit (CN-unit).} 
\label{fig30}
}
\end{figure}

In this section, notice the tilings shown in Section~\ref{section3} that contain CN-units 
(see Figure~\ref{fig30}). As described in \cite{S_and_A_2017}, CN-units in 
the tiling can be freely reversed. With this operation, they are a new convex 
pentagon tiling (they will be new heptiamond tilings).

\begin{figure}[htbp]
 \centering\includegraphics[width=14.5cm,clip]{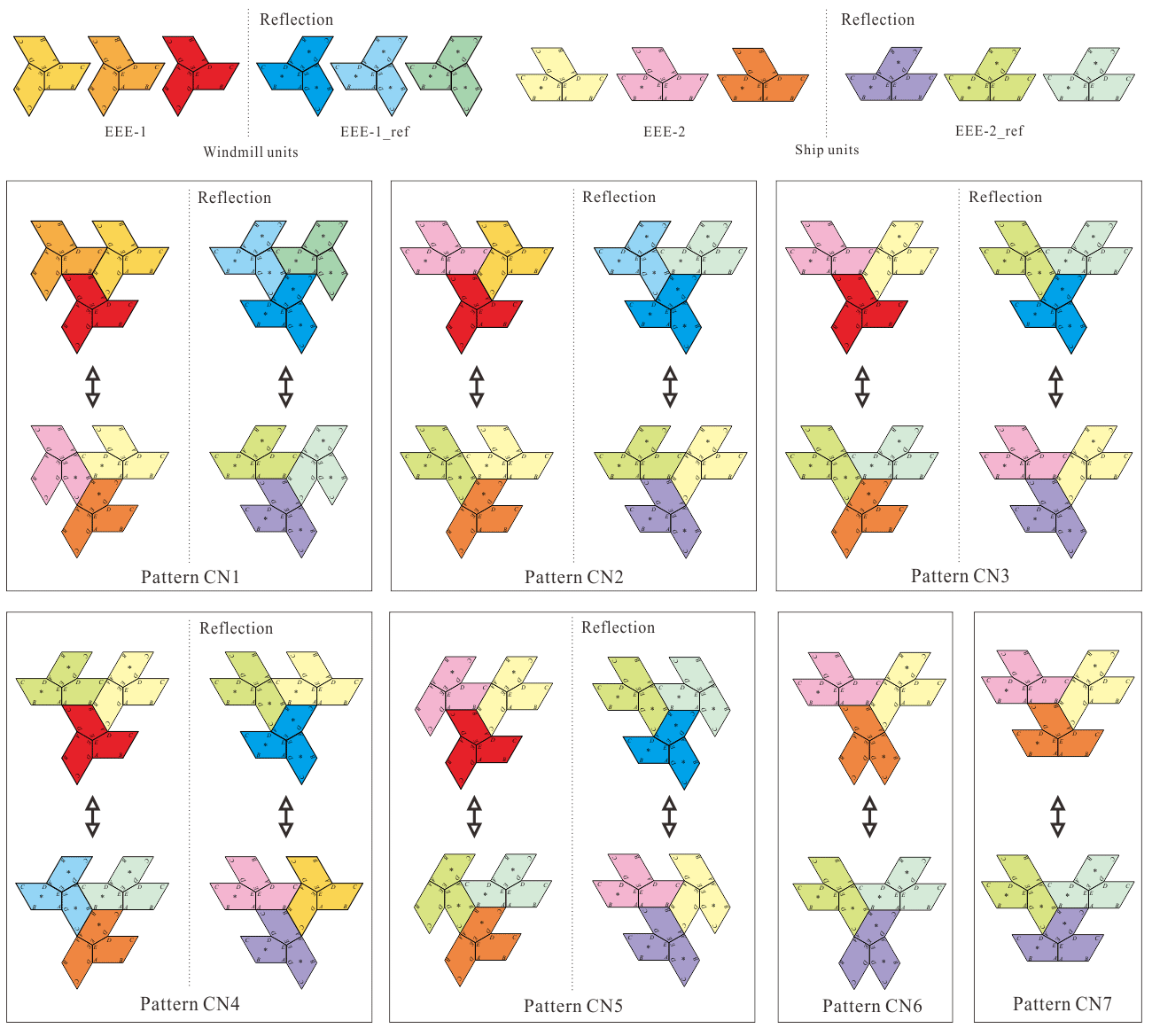} 
  \caption{{\small 
Seven unique patterns of CN-units that can be formed using 
windmill units and ship unit.} 
\label{fig31}
}
\end{figure}

For a CN-unit that can be formed using windmill units and ship units, there 
are seven unique patterns of pairs for ACN-units (anterior CN-unit) and 
PCN-units (posterior CN-unit) shown in Figure~\ref{fig31}. The Patterns 
that used in Rice1995-tiling and Figure 54 in \cite{S_and_A_2017} correspond 
to Pattern CN1 and Pattern CN2 in Figure~\ref{fig31}, respectively. For Patterns 
CN6 and CN7, even if the CN-unit is reversed, they are always formed 
by only the ship unit (however, the number of anterior and posterior 
convex pentagons is changed).

\begin{figure}[htbp]
 \centering\includegraphics[width=15cm,clip]{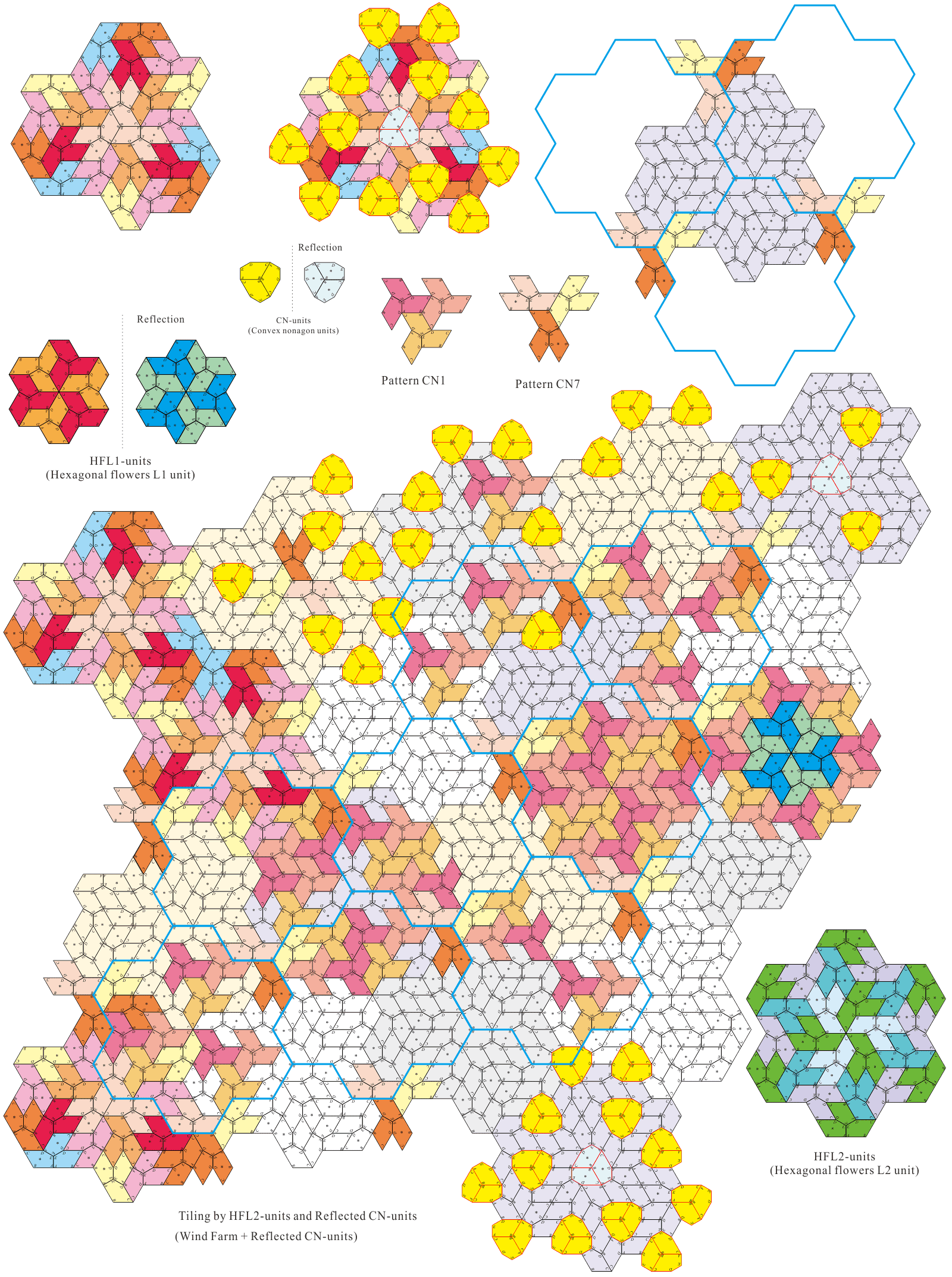} 
  \caption{{\small 
Tiling of Figure~\ref{fig23} with reversed CN-units (Tiling by 
HFL2-units).} 
\label{fig32}
}
\end{figure}

\begin{figure}[htbp]
 \centering\includegraphics[width=15cm,clip]{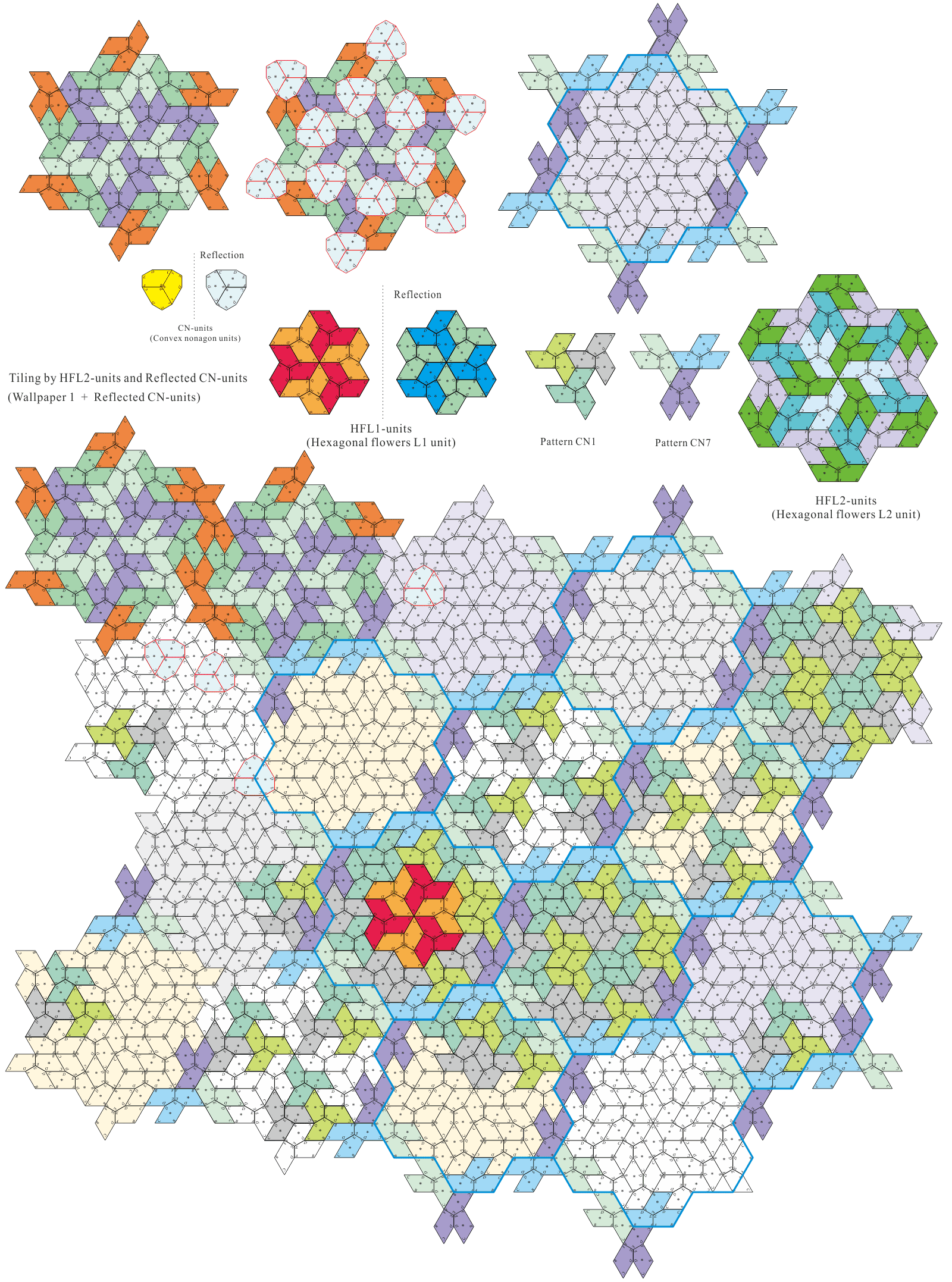} 
  \caption{{\small 
Tiling of Figure~\ref{fig24} with reversed CN-units (Tiling by 
HFL2-units).} 
\label{fig33}
}
\end{figure}

\begin{figure}[htbp]
 \centering\includegraphics[width=15cm,clip]{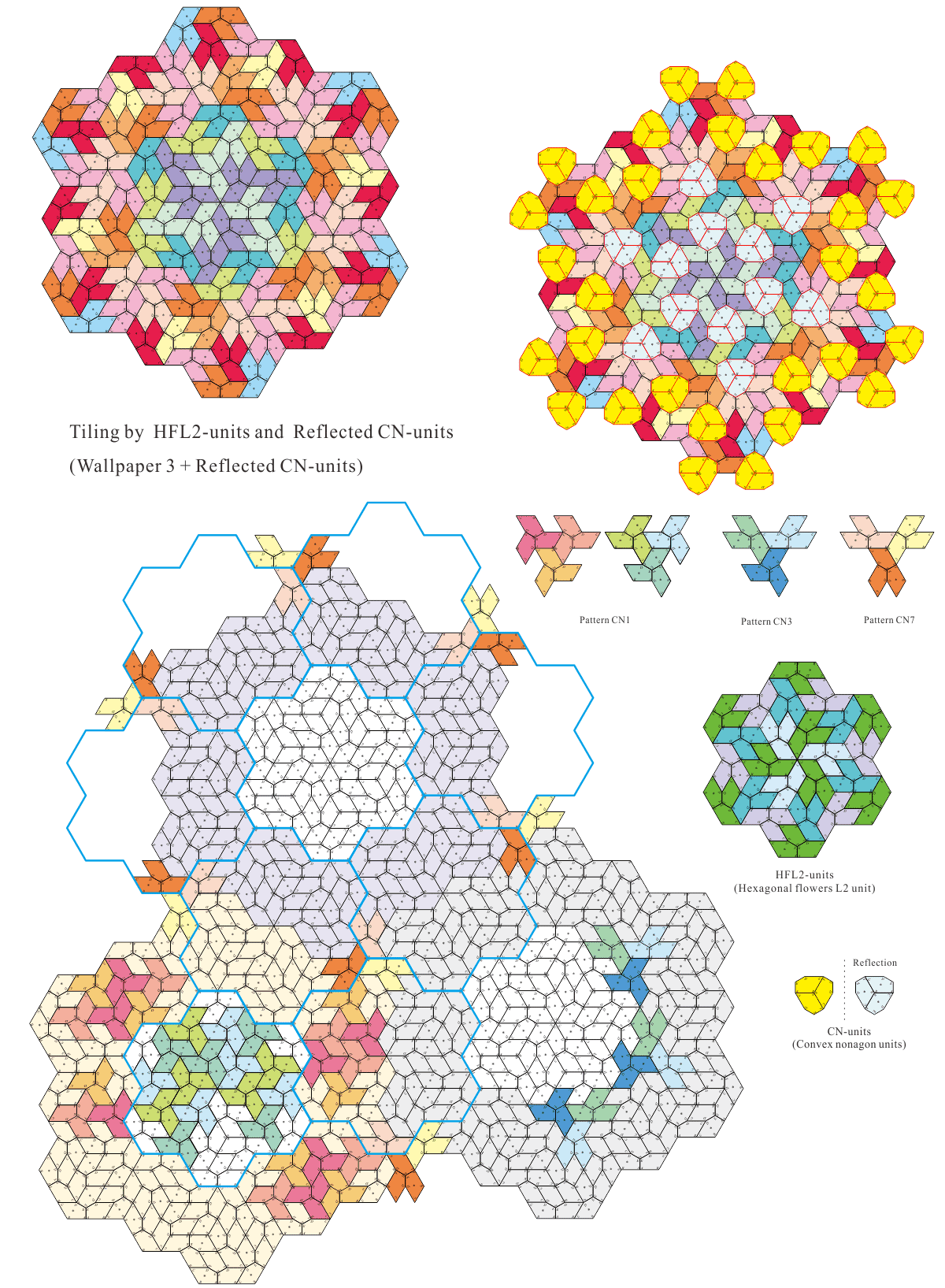} 
  \caption{{\small 
Tiling of Figure~\ref{fig26} with reversed CN-units (Tiling by 
HFL2-units).} 
\label{fig34}
}
\end{figure}

First, the tilings of Figures~\ref{fig23}, \ref{fig24}, and \ref{fig26} (Wind Farm, 
Wallpaper 1, and Wallpaper 3) where the explanation in the previous section 
was passed to this section are shown. These are tilings of HFL2-units 
formed with only the ship unit and can be formed by reversing CN-units on 
the boundaries of HFL2-units. Figures~\ref{fig32}, \ref{fig33}, and \ref{fig34}
 show the properties. As shown in \cite{S_and_A_2017} (see 
Figure~\ref{fig41}), the HFL2-units can also be formed by windmill units and ship 
units, so that more various tiling with convex pentagon can be generated by 
using them. In addition, CN-units exist in Rice1995-tiling, tilings with 
HFL2-units, or cases that connected those two tilings by only the ship units 
in Figures~\ref{fig9}, \ref{fig28}, and \ref{fig29}. Therefore, these tilings can be 
shifted to tilings with reversed CN-units. The convex pentagon tilings in the 
above case are variations of tilings which can be created according to these 
considerations of \cite{S_and_A_2017}. 

\begin{figure}[htbp]
 \centering\includegraphics[width=15cm,clip]{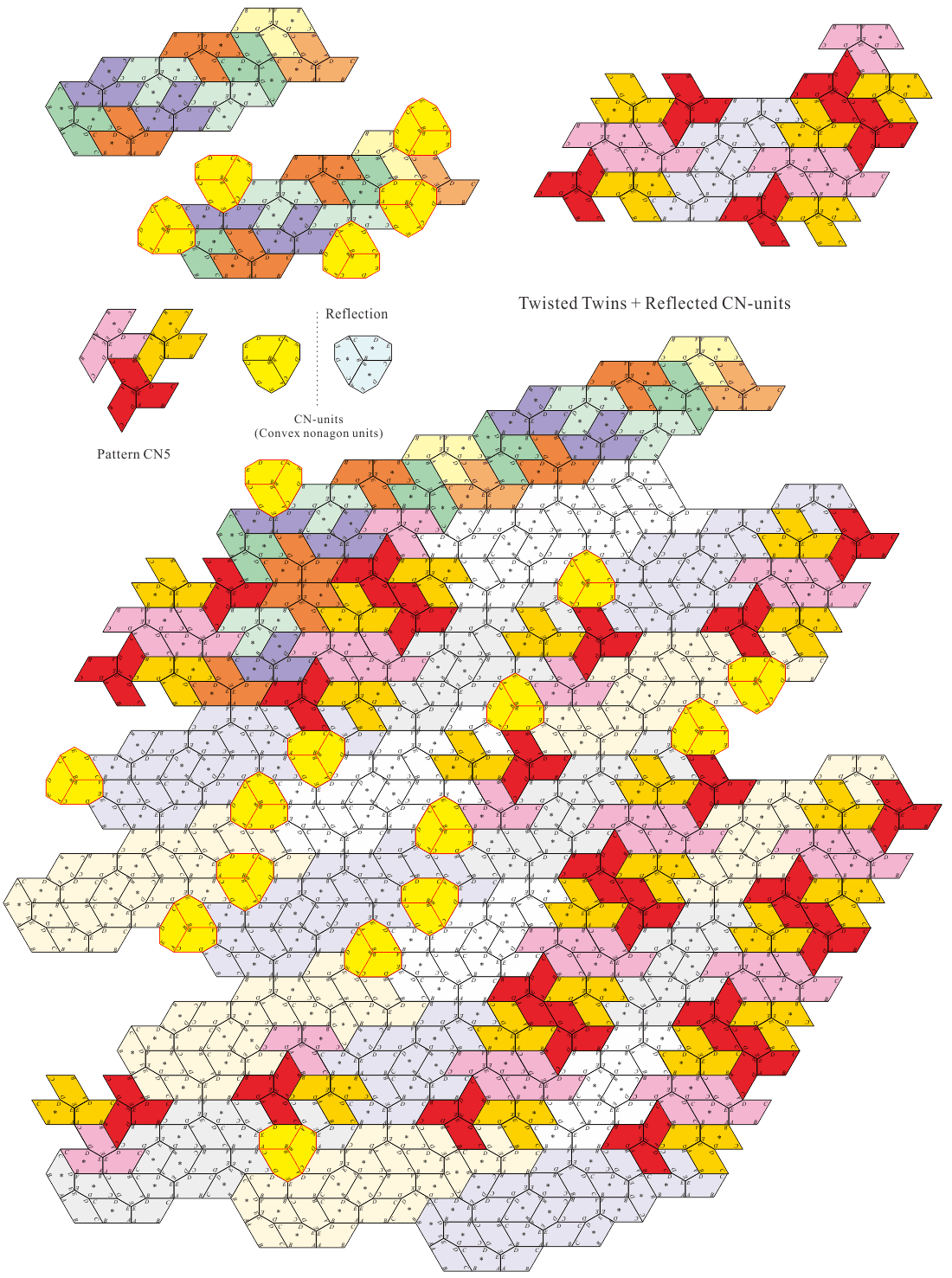} 
  \caption{{\small 
Tiling of Figure~\ref{fig14} with reversed CN-units (New convex pentagon 
tiling).} 
\label{fig35}
}
\end{figure}

\begin{figure}[htbp]
 \centering\includegraphics[width=15cm,clip]{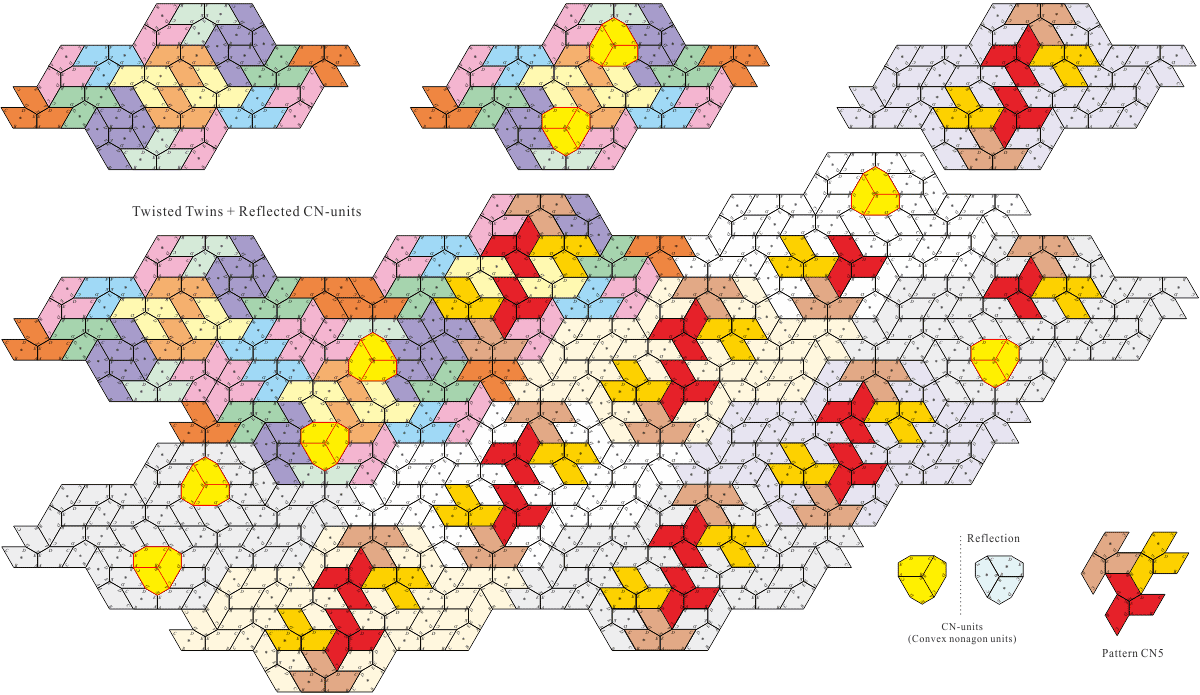} 
  \caption{{\small 
Tiling of Figure~\ref{fig17} with reversed CN-units (New convex pentagon 
tiling).} 
\label{fig36}
}
\end{figure}

\begin{figure}[htbp]
 \centering\includegraphics[width=15cm,clip]{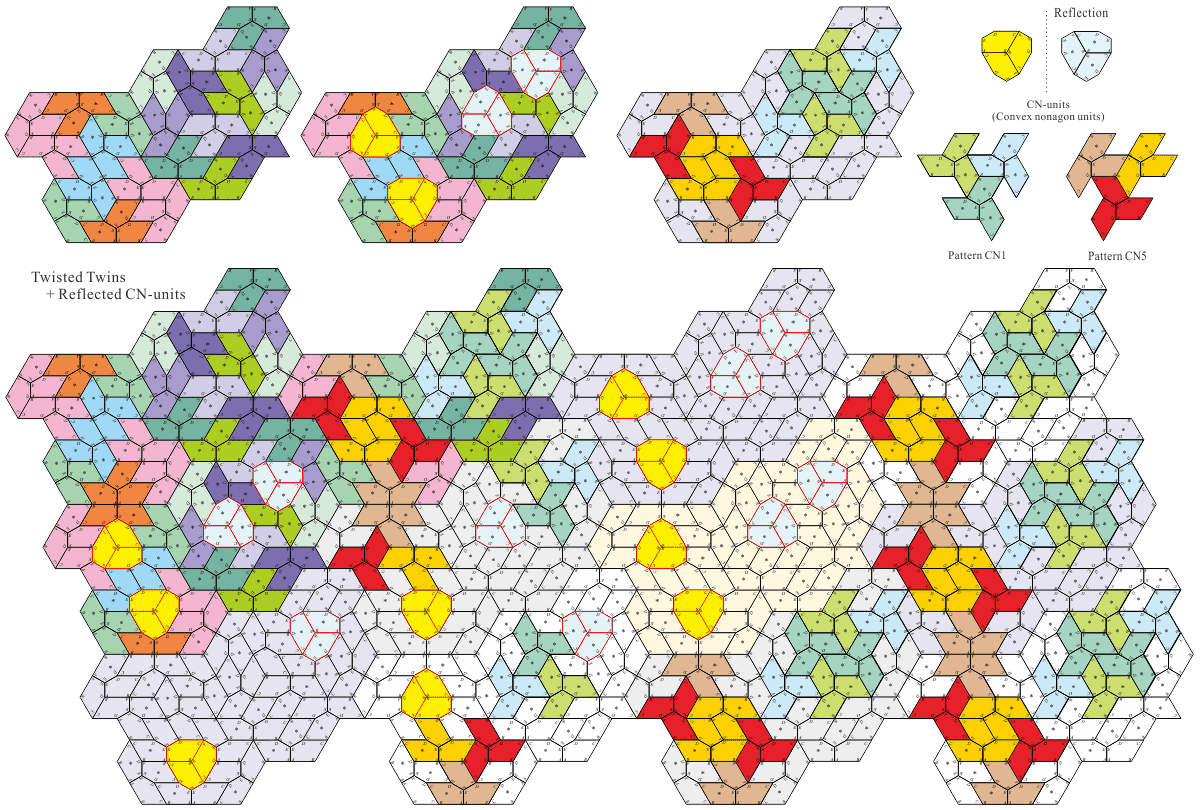} 
  \caption{{\small 
Tiling of Figure~\ref{fig18} with reversed CN-units (New convex pentagon 
tiling).} 
\label{fig37}
}
\end{figure}

\begin{figure}[htbp]
 \centering\includegraphics[width=14cm,clip]{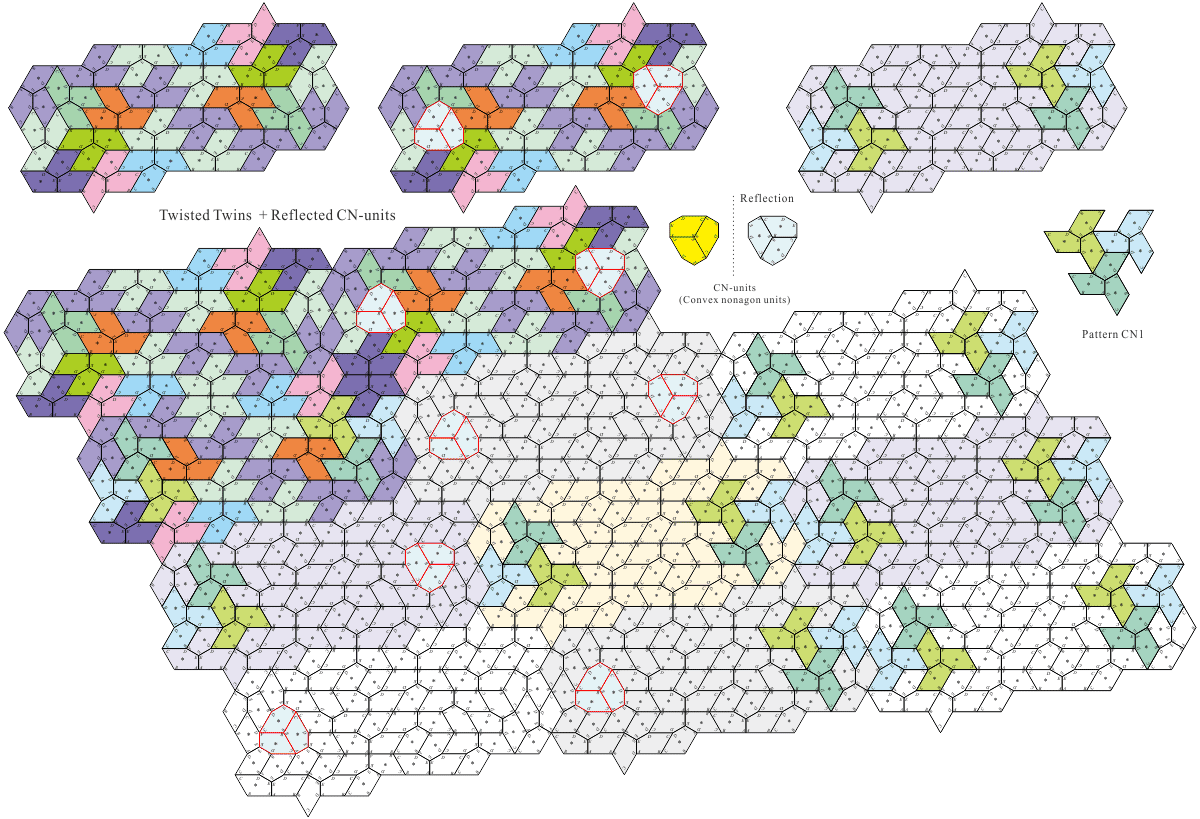} 
  \caption{{\small 
Tiling of Figure~\ref{fig19} with reversed CN-units (New convex pentagon 
tiling).} 
\label{fig38}
}
\end{figure}

\begin{figure}[htbp]
 \centering\includegraphics[width=14cm,clip]{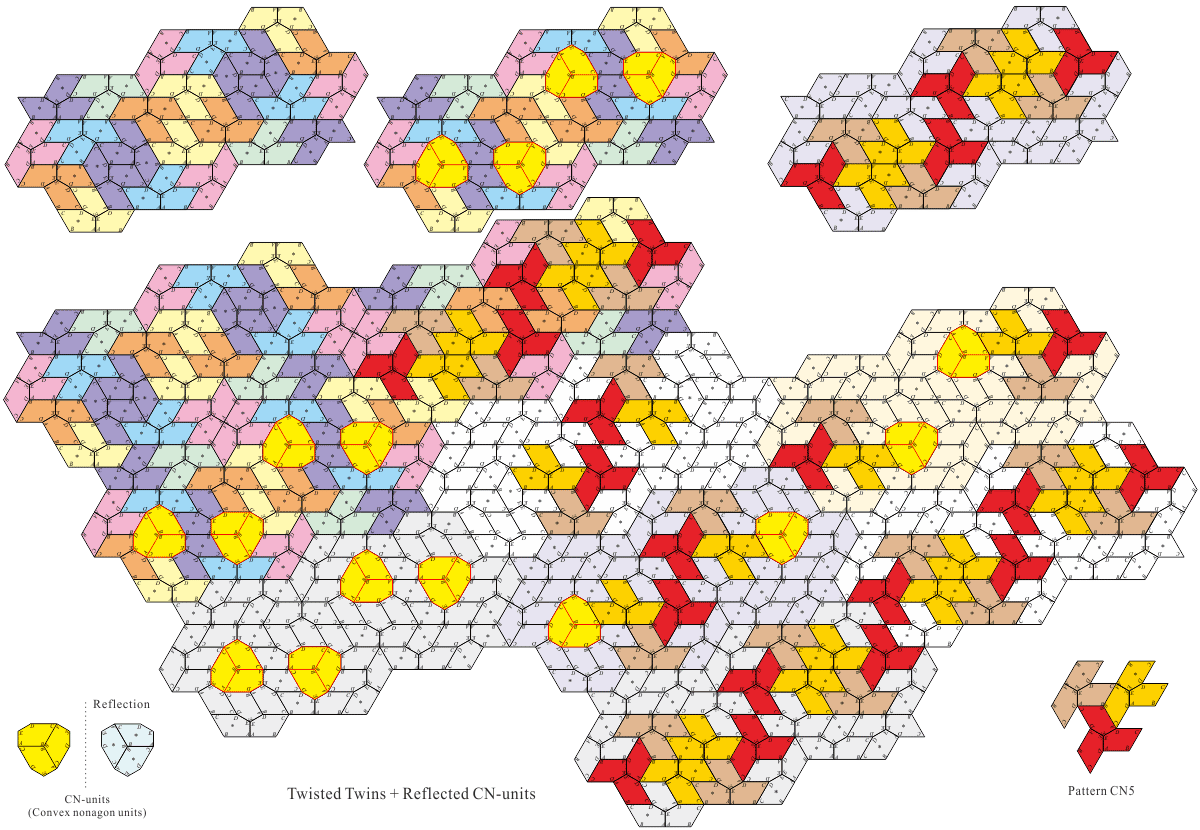} 
  \caption{{\small 
 Tiling of Figure~\ref{fig20} with reversed CN-units (New convex pentagon 
tiling).} 
\label{fig39}
}
\end{figure}

\begin{figure}[htbp]
 \centering\includegraphics[width=15cm,clip]{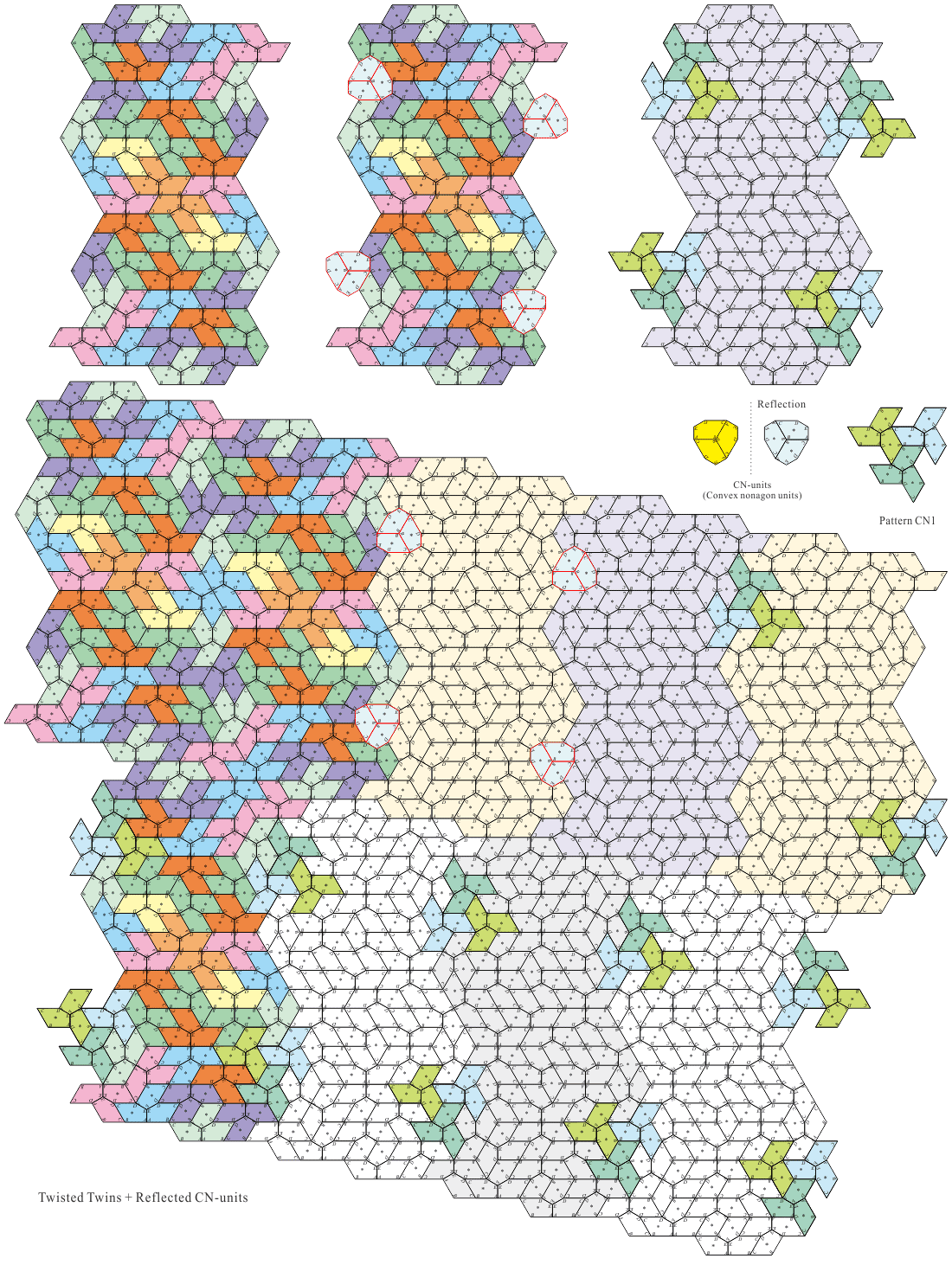} 
  \caption{{\small 
Tiling of Figure~\ref{fig21} with reversed CN-units (New convex pentagon 
tiling).} 
\label{fig40}
}
\end{figure}

\begin{figure}[htbp]
 \centering\includegraphics[width=14.5cm,clip]{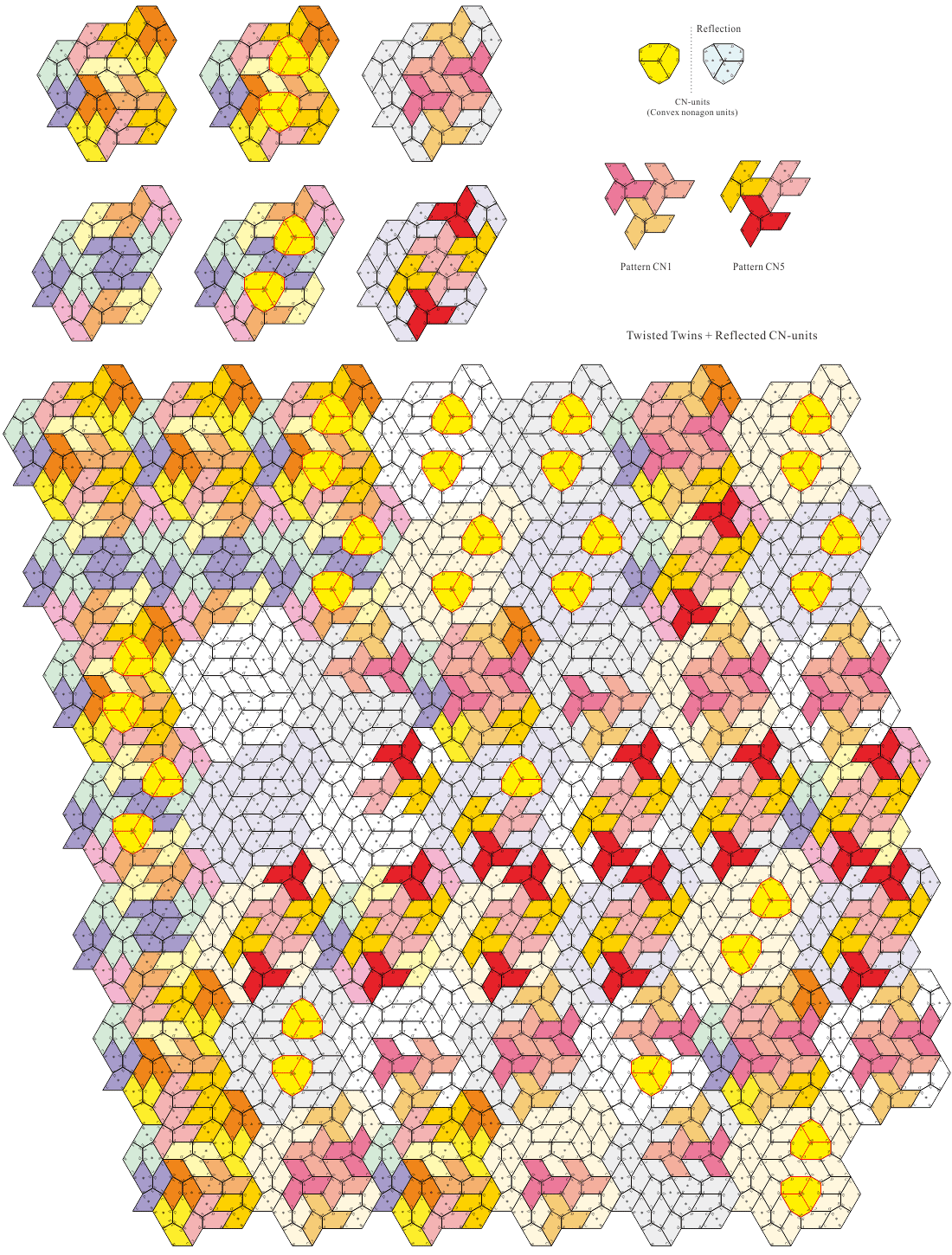} 
  \caption{{\small 
Tiling of Figure~\ref{fig22} with reversed CN-units (New convex pentagon 
tiling).} 
\label{fig41}
}
\end{figure}

\begin{figure}[htbp]
 \centering\includegraphics[width=15cm,clip]{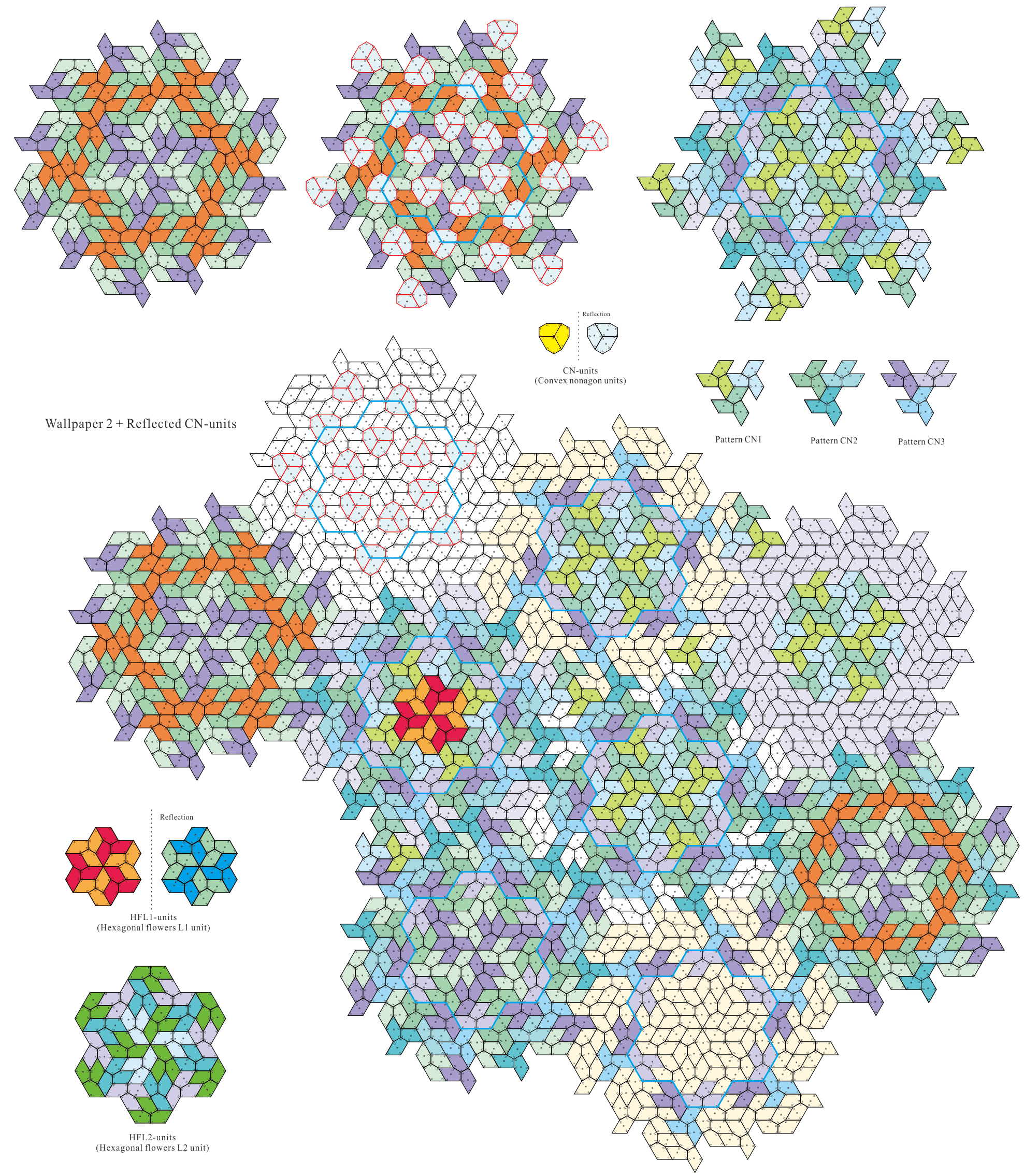} 
  \caption{{\small 
Tiling of Figure~\ref{fig25} with reversed CN-units (New convex pentagon 
tiling).} 
\label{fig42}
}
\end{figure}

\begin{figure}[htbp]
 \centering\includegraphics[width=15cm,clip]{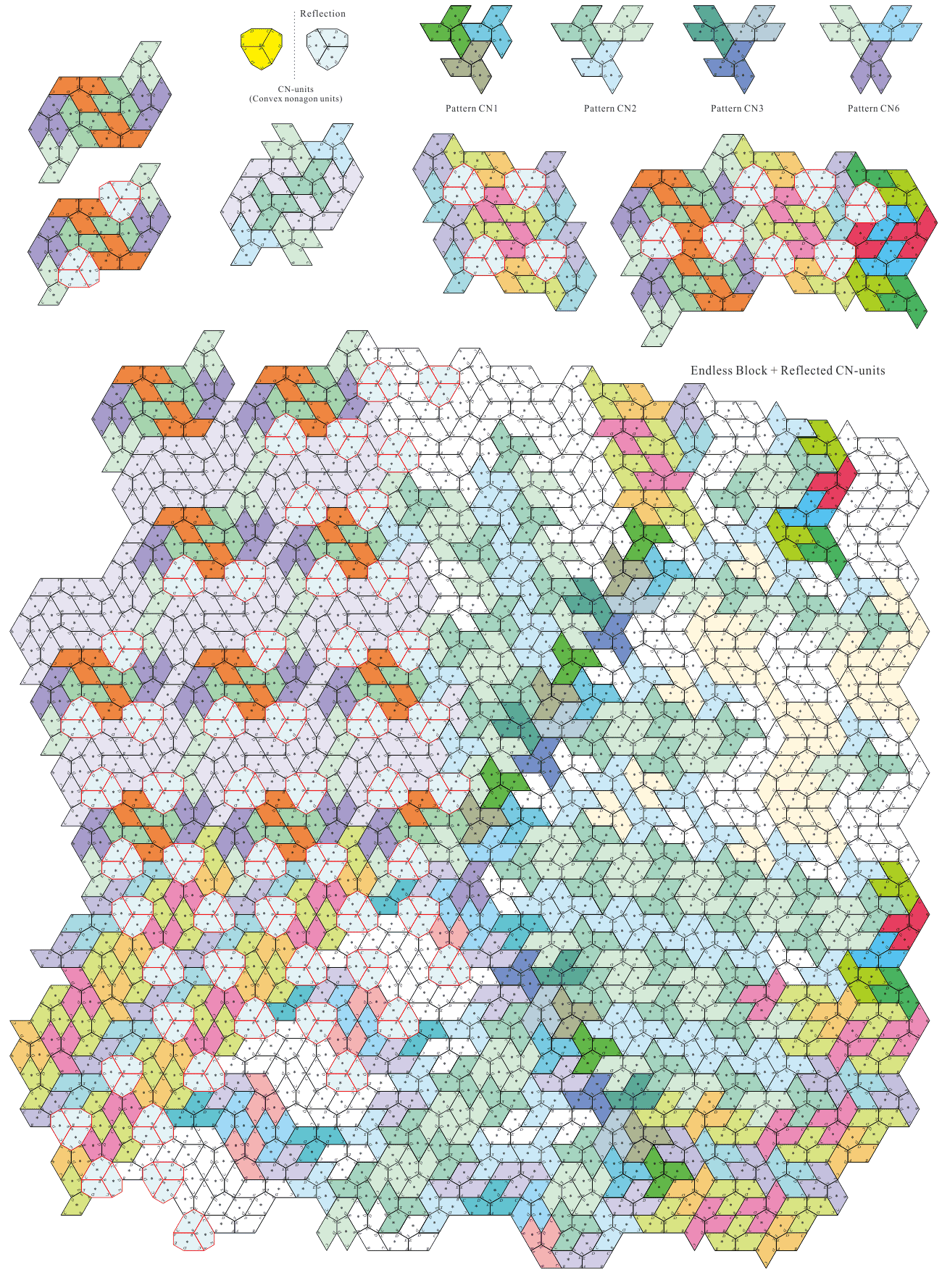} 
  \caption{{\small 
Tiling of Figure~\ref{fig27} with reversed CN-units (New convex pentagon 
tiling).} 
\label{fig43}
}
\end{figure}

On the other hand, the tilings of Figures~\ref{fig14}, \ref{fig17}-\ref{fig22}, 
and \ref{fig25} have CN-units and are new convex pentagon tilings. 
That is, reversing CN-units in these tilings will result in new convex pentagon 
tilings. Because CN-units in tilings can be reversed freely, there will be infinite 
periodic and non-periodic tilings. Therefore, since it is impossible to draw all tilings, 
in this manuscript, examples of each tiling with reversed CN-units are shown 
(see Figures~\ref{fig35}-\ref{fig43}).

\renewcommand{\figurename}{{\small Figure.}}

\renewcommand{\figurename}{{\small Figure.}}

\renewcommand{\figurename}{{\small Figure.}}

\renewcommand{\figurename}{{\small Figure.}}

\renewcommand{\figurename}{{\small Figure.}}

\renewcommand{\figurename}{{\small Figure.}}

\renewcommand{\figurename}{{\small Figure.}}

\renewcommand{\figurename}{{\small Figure.}}

\renewcommand{\figurename}{{\small Figure.}}

\renewcommand{\figurename}{{\small Figure.}}

\renewcommand{\figurename}{{\small Figure.}}

\renewcommand{\figurename}{{\small Figure.}}

\renewcommand{\figurename}{{\small Figure.}}

\renewcommand{\figurename}{{\small Figure.}}

\section{Conclusion }

Hindriks' search results are very interesting and his content is 
substantial. However, since Hindriks considers only heptiamonds, there are 
no contents of convex pentagon tilings as presented here in Hindriks' site. 
There are no contents of properties like reversing CN-units. Therefore, 
there is no such content that the tilings will change by reversing CN-units.

Neither the authors nor Hindriks can examine everything. Therefore, there is 
a sufficient possibility that new tiling (a new fundamental region) will be 
found in the future\footnote{ In response to our results, Johannes Hindriks 
discovered new heptiamond tilings with windmill units and ship units. 
We converted those heptiamond tilings into convex pentagon tilings. 
As for the results, refer to Sugimoto's site:  
\url{http://tilingpackingcovering.web.fc2.com/th-pentagon/h-201703-e.html} and
\url{http://tilingpackingcovering.web.fc2.com/th-pentagon/h-201704-e.html}}.

\end{document}